\documentclass[preprint,12pt]{elsarticle2}
\usepackage[english]{babel}
\usepackage[T1]{fontenc}
\usepackage{amsmath,amsthm,amssymb}
\usepackage{amsfonts}
\usepackage{appendix}
\usepackage{geometry}
\usepackage{bbold}
\usepackage{mathrsfs}
\usepackage{graphicx}
\usepackage{tikz}
\usepackage{array}
\usepackage{natbib}
\usepackage{makecell}
\usepackage{hyperref}

\graphicspath{{figures}{figures/1D}{figures/2D}{figures/non_local}}
\usetikzlibrary{arrows}
\usetikzlibrary{calc}
\usetikzlibrary {arrows.meta} 

\journal{Journal of Computational Physics}
\usetikzlibrary{shapes.arrows,chains,patterns,matrix,arrows}
\newcommand*{\crochet}[1]{\left<#1\right>}
\newcommand*{\sphere}[0]{\mathbb{S}^2}
\newcommand*{\ccrochet}[1]{\left<\left<#1\right>\right>}

\newcommand*{\indi}[3]{#1 \mathbb{1}_{#3 > 0} + #2 \mathbb{1}_{#3 < 0}}

\newcommand*{\derivee}[2]{\dfrac{\mathrm{d}#1}{\mathrm{d}#2}}

\newcommand*{\mf}[1]{\mathbf{#1}}

\DeclareMathOperator{\sign}{sign}

\newtheorem{prop}{Proposition}[section]

\newtheorem{defi}{Definition}[section]

\begin{document}
\begin{frontmatter}
\title{An asymptotic preserving kinetic scheme for the M1 model of non-local thermal transport for two-dimensional structured and unstructured meshes}
\author[1]{Jean-Luc Feugeas}
\ead{jean-luc.feugeas@u-bordeaux.fr}

\author[2]{Julien Mathiaud}
\ead{julien.mathiaud@univ-rennes.fr}

\author[3]{Luc Mieussens}
\ead{Luc.Mieussens@math.u-bordeaux.fr}

\author[1]{Thomas Vigier\corref{cor1}}
\ead{thomas.vigier@u-bordeaux.fr}

\affiliation[1]{organization={CELIA, University of Bordeaux, CNRS, CEA},
            addressline={UMR 5107}, 
            city={Talence},
            postcode={F-33405},
            country={France}}
\affiliation[2]{organization={University of Rennes, CNRS, IRMAR},
            addressline={UMR 6625}, 
            city={Rennes},
            postcode={F-35700},
            country={France}}
\affiliation[3]{organization={University of Bordeaux, Bordeaux INP, CNRS, IMB},
            addressline={UMR 5251}, 
            city={Talence},
            postcode={F-33400},
            country={France}}
\cortext[cor1]{Corresponding author}

\begin{abstract}
The M1 moment model for electronic transport is commonly used to describe non-local thermal transport effects in laser-plasma simulations. In this article, we propose a new asymptotic-preserving scheme based on the Unified Gas Kinetic Scheme (UGKS) for this model in two-dimensional space. This finite volume kinetic scheme follows the same approach as in our previous article \cite{feugeas2023asymptotic} and relies on a moment closure, at the numerical scale, of the microscopic flux of UGKS. The method is developed for both structured and unstructured meshes, and several techniques are introduced to ensure accurate fluxes in the diffusion limit. A second-order extension is also proposed. Several test cases validate the different aspects of the scheme and demonstrate its efficiency in multiscale simulations. In particular, the results demonstrate that this method accurately captures non-local thermal effects.
\end{abstract}
\end{frontmatter}

\section{Introduction}
In the context of plasma physics and in particular in laser plasma interactions, accurately describing the electron heat flux is crucial. Indeed, this flux needs to be expressed in terms of the macroscopic quantities (density, velocity, pressure and temperature) that govern the plasma evolution at the macroscopic scale. Close to thermodynamic equilibrium, the Spitzer-Harm model \cite{spitzer1953transport} provides a suitable closure. However, in most practical situations such as inertial confinement fusion and radiotherapy, this description is insufficient as the presence of strong temperature gradients induces non-local kinetic effects that cannot be recovered with this law. As a matter of fact, non-local models have been developed to correctly describe electron heat fluxes (for example see \cite{schurtz2000nonlocal,nicolai2006practical,chrisment2022analysis}). These models rely on a kinetic description of the electron population to evaluate this flux. This approach significantly increases the computational complexity. Nevertheless, an intermediate description at the mesoscopic scale seems to be sufficient in practice to reproduce most kinetic effects. In particular directional moment models, such as the M1 one, are commonly used \cite{del2017benchmark,chrisment2023deterministic}. 

In general, moment models rely on a moment hierarchy of the kinetic equation to reduce the number of kinetic variables and hence the problem dimension. The integration process leads to a last unknown flux that requires closure. Given a specific ansatz for the distribution function shape, the hierarchy can be closed and this flux can be expressed in terms of the model unknowns. The chosen ansatz entirely determines the kinetic effects that can be recovered and is usually chosen based on physical argument. Under the small anisotropy hypothesis, the Pn moment model relies on a truncated expansion of the distribution function in spherical harmonics. The associated closure is linear but does not ensure the underlying distribution function positivity and hence the moments realizability \cite{modelcoll}. In contrast, the Mn moment model relies on an entropic argument to capture the ansatz. The chosen distribution function minimizes an entropy functional (usually the system entropy which is in this case the Boltzmann one) under constraints on the moments. This procedure allows to recover a positive distribution function and to define a hyperbolic model which ensures the moments realizability \cite{levermore1996moment,dubroca1999etude,alldredge2012high}.  

From a numerical point of view, developing accurate numerical schemes for such models of electron transport is challenging as they must capture the smallest microscopic scale in the domain, which constrains the space discretization and thus the time steps for stability reasons. The Knudsen number (denoted by $\epsilon$) is defined as the ratio between the mean free path of the particles and a macroscopic length and is representative of the collision regime. As this dimensionless number tends to zero, a global diffusion mechanism emerges at the macroscopic scale. Numerical schemes for the moment models are not necessarily consistent with the macroscopic model in this limit. To address this issue, asymptotic-preserving (AP) numerical schemes have been introduced in the literature. These schemes are specially designed to be consistent with the limit model while being uniformly stable with $\epsilon$ (see \cite{larsen1987,larsen1989,jin1991,jin1993,jin2000,klar1998,buet2002diffusion,klar2001numerical,lemou2008new,bennoune2008uniformly,carrillo2008simulation,carrillo2008numerical,gosse2011transient,lafitte2012asymptotic} for examples of AP schemes in a general context). In particular, with regard to the M1 model, AP schemes have been developed, for instance in \cite{guisset2018admissible,turpault,berthon2007asymptotic,feugeas2024asymptotic}.

The Unified Gas Kinetic Scheme (UGKS) is an AP scheme originally developed by Xu and Huang in 2010 in the context of rarefied gas dynamics \cite{kxu2010ugks}. Since then, it has been further improved and the general ideas have been applied to complex gas flows \cite{liu2017unified} (see \cite{zhu2021first} for other references). The UGKS was also extended to linear models with the diffusion limit in \cite{mieussens2013,sun2015radiativetransfer}. 

This paper follows up on our previous study on the linear transport model, where the UGKS was used to develop an AP numerical scheme for the M1 model of linear transport \cite{feugeas2024asymptotic}. Despite the rather simple context, this work proposes a general procedure to obtain AP schemes for moment models. The goal of this paper is to apply the same approach to the more relevant context of the non-local thermal transport theory. This new context introduces two difficulties; the scheme is developed in two-dimensional space and the collision operator acts towards a non-trivial equilibrium. In spite of this, we follow here the same generic approach to obtain a numerical scheme for the M1 moment model of the electron transport equation. The central idea is to apply the M1 closure at the numerical level in the numerical approximation (given by the UGKS) of the kinetic equation. We prove that this procedure inherits from the asymptotic-preserving property of the UGKS and hence captures the correct diffusion regime. Moreover, we suggest a second order extension that does not compromise the AP property. Additionally, we show that different diffusion schemes can be easily recovered. This work is first presented on structured meshes and then extended on unstructured meshes. In the latter case, particular attention is given in the UGKS construction to obtain proper diffusion schemes in the limit. 

The outline of our paper is as follows. First in section \ref{sec:model}, we introduce the kinetic model and construct the associated entropic moment model. The main properties of both models are also highlighted. Next in section \ref{sec:ugks}, the UGKS is constructed for this specific kinetic equation and a new scheme for the M1 moment model is derived from it (several extensions are proposed and some numerical difficulties are discussed). In section \ref{sec:uns}, the UGKS and its extension to the M1 model are developed for unstructured meshes. A method is proposed to obtain the diamond scheme in the diffusion limit. Finally, the schemes are evaluated on several test cases in section \ref{sec:results}.

\section{The entropic closure for the electron transport equation \label{sec:model}}
    \subsection{The electron transport equation}

The electron transport equation is a kinetic equation that describes the evolution of the particle distribution function $f$ as a function of time $t$, of space position $\mathbf{x}$ in $\mathscr{D}$ an open set of $\mathbb{R}^3$, of speed $v$ in $\mathbb{R}_+$ and of velocity direction $\boldsymbol\Omega$ in $\sphere$ the unit sphere in three-dimensional space:

\begin{equation} \label{equation_base_non_adim}
    \partial_t f + v \boldsymbol{\Omega}\cdot \boldsymbol{\nabla_\mathbf{x}} f = \dfrac{1}{\tau}\left( M[f]-f\right).
\end{equation}
The electron distribution $f(t,\mathbf{x},v,\boldsymbol\Omega)$ represents the density of electrons in the space volume $v^2\mathrm{d}v\mathrm{d}\boldsymbol\Omega$ at a certain time $t$. From a physical point of view, the kinetic equation (\ref{equation_base_non_adim}) describes the evolution in time of the electron distribution in the absence of external forces. Thus, the influence of electromagnetic fields through the Lorentz force is neglected. The left-hand side of the equation is total derivative in time of the distribution function and expresses the electron advection at speed $v$ in the direction $\mathbf \Omega$. The right hand side accounts for the particle collisions. Those interactions are modelled with the Bhatnagar-Gross-Krook (BGK) operator and represents the rate of change of the distribution function. It acts as a relaxation term towards the equilibrium state, which is the Maxwellian distribution $M[f]$, within a characteristic collision time $\tau$. This operator is an approximation of the non-linear Boltzmann one, while still preserving some fundamental properties such as the relaxation to the correct equilibrium, the mass and energy conservation, as well as the entropy dissipation. 

To study the diffusion regime and for computational reasons, it is convenient to work with the non-dimensional equation. Several characteristic quantities are introduced: $t^*$ a characteristic time, $L$ a characteristic macroscopic length, $c$ a characteristic speed and $\tau^*$ a characteristic collision time. The non-dimensional equation is
\begin{equation} \label{equation_base}
   \partial_t f + \dfrac{v}{\eta} \boldsymbol{\Omega}\cdot \boldsymbol{\nabla}_\mathbf{x} f = \nu\left(M[f]-f\right),
\end{equation}
where $\epsilon$ is the Knudsen number defined as the ratio between the mean free path ($\lambda=c\tau^*$) and the macroscopic length $L$ and $\eta$ is an analogue of the Strouhal number and is defined as the ratio between the macroscopic velocity and the microscopic one:
\begin{equation} \label{admi_var}
    \epsilon=\dfrac{\lambda}{L}, \quad \eta=\dfrac{L/\tau^*}{c}.
\end{equation}
In equation (\ref{equation_base}), the non-dimension Maxwellian $M[f]$ is the isotropic distribution of equilibrium defined as
\begin{equation} \label{max}
    M[f](v)=\dfrac{\rho}{(2\pi T)^{3/2}} \exp \left(-\dfrac{v^2}{2T}\right),
\end{equation}
and the collision frequency is defined as
\begin{equation}
    \nu=\dfrac{C}{\epsilon\eta}\dfrac{\rho}{T^{3/2}},
\end{equation}
where $C$ is a physical constant, $\rho$ is the macroscopic number density and $T$ is temperature. The energy density is $q=\frac{3}{2}\rho T$. Both the Maxwellian and the collision frequency only depend on the macroscopic state $\mf{W}=\begin{pmatrix}\rho & q\end{pmatrix}^T$ of the system. Those macroscopic variables are moments of the distribution function:
\begin{equation}
    \mathbf{W}=\ccrochet{\boldsymbol\Psi f}=\int_{\sphere}\int_{\mathbb{R}_+}\boldsymbol\Psi(v) f(\cdot,\cdot,v,\boldsymbol\Omega)v^2\mathrm{d}v\mathrm{d}\Omega ,
\end{equation}
where $\boldsymbol\Psi(v)=\begin{pmatrix}1&\frac{1}{2}v^2\end{pmatrix}^T$ is the vector of collisional invariants. Those macroscopic quantities are conserved as
\begin{equation}
    \ccrochet{\boldsymbol\Psi(M[f]-f)}=0.
\end{equation}
Thus, two conservation laws can be written:
\begin{equation}\label{conservation_macro}
\partial_t \mf{W} + \ccrochet{ \dfrac{v}{\eta}\boldsymbol\Psi \left(\boldsymbol\Omega \cdot \boldsymbol\nabla_\mf{x} f\right)}=0.
\end{equation}
In the following, it is convenient to introduce the opacity $\sigma$ which is defined so that
\begin{equation}
    \nu = \dfrac{\sigma}{\epsilon\eta}.
\end{equation}

\subsection{Asymptotic regimes}
As the Knudsen number tends to zero, the collision mechanism predominates in comparison with particle transport. As a consequence, the regime is no longer kinetic and the distribution functions tend to be Maxwellians. At the diffusion scale (when $\eta=\epsilon$), a global diffusion mechanism emerges on the macroscopic variables. A Hilbert expansion of the distribution function and of the Maxwellian allows to show that the macroscopic variables $\mf{W}$ satisfy a coupled system of diffusion at the first order in $\epsilon$:
\begin{equation} \label{diffusion}
    \begin{cases}
        \partial_t \rho &= \boldsymbol\nabla_\mathbf{x} \cdot \left(\dfrac{2}{3\sigma}\boldsymbol\nabla_\mathbf{x} q\right) + \mathcal{O}(\epsilon)\\
        \partial_t q &= \boldsymbol\nabla_\mathbf{x} \cdot \left(\dfrac{5}{3}\dfrac{2}{3\sigma}\boldsymbol\nabla_\mathbf{x} \left(\dfrac{q^2}{\rho}\right)\right) + \mathcal{O}(\epsilon)\\
    \end{cases},
\end{equation}
The mathematical study of this system is rather delicate. However it can be shown that this system still ensures the entropy dissipation property.

Conversely, in the free transport regime ($\eta$ constant and $\epsilon\to +\infty$), the electrons are advected at velocity $\frac{v}{\eta}$ in direction $\boldsymbol\Omega$ without interacting with each other. Hence, as $\epsilon$ tends to infinity, the usual linear advection equation without is a source term is recovered:
\begin{equation} \label{equation_transp}
   \partial_t f + \dfrac{v}{\eta} \boldsymbol{\Omega}\cdot \boldsymbol{\nabla}_\mathbf{x} f = 0.
\end{equation}

    \subsection{The M1 moment closure}
Solving equation (\ref{equation_base}) is expensive due to the dimension of the problem. However, in most physical applications, assumptions can be made on the shape of the distribution function. For instance, close to equilibrium state, the distribution is nearly isotropic (it only depends on the speed $v$). On the other hand, in some non-equilibrium situations (such as beams), all particles tends to have the same direction. Assuming the particle distribution in direction, it is possible to develop reduced models in $\boldsymbol\Omega$. These models capture some kinetic effects while significantly lowering the computational cost by lowering the problem's dimensionality. \\

This first step to elaborate such a model is to get ride of this kinetic variable by forming an angular moment hierarchy of the kinetic equation. The simplest hierarchy which enables the restoration of an angular anisotropy is obtained by integrating the kinetic equation (\ref{equation_base}) against the vector $\mathbf{m}(\boldsymbol\Omega)=\begin{pmatrix} 1 & \boldsymbol\Omega\end{pmatrix}^T \in \mathbb{R}^4$ with respect to the velocity direction $\boldsymbol\Omega$. First, the distribution function angular moments are
\begin{subequations}
    \begin{align}
        f_0  &=\crochet{ f(\cdot,\boldsymbol\Omega)}=\int_{\sphere}f(\cdot,\boldsymbol\Omega) \mathrm{d}\Omega \\
        \mf{f}_1&=\crochet{\boldsymbol\Omega f(\cdot,\boldsymbol\Omega)} \\
        \mf{f}_2&=\crochet{\boldsymbol\Omega \otimes \boldsymbol\Omega f(\cdot,\boldsymbol\Omega)}
    \end{align}
\end{subequations}
where the tensor $\mathbf{f}_i$ has rank $i$. In the following, integrals of any function of $\boldsymbol\Omega$ over the sphere are noted $\crochet{\cdot}$. Finally, the moment hierarchy is
\begin{equation}\label{model_moment}
    \partial_t\mathbf{U} + \partial_x \mathbf{F}(\mathbf{U}) = \nu \mathbf{S}(\mathbf{U}),
\end{equation}
where $\mathbf{U}=\crochet{\mf{m}f}=\begin{pmatrix} f_0 & \mathbf{f}_1 \end{pmatrix}^T \in \mathbb{R}^4$ is the vector of conservative variables, $\mathbf{F}=\frac{v}{\eta} \begin{pmatrix} \mathbf{f}_1 & \mathbf{f}_2 \end{pmatrix}^T$ is the flux vector and $\mathbf{S}(\mathbf{U})=\begin{pmatrix} M_0[\mathbf{U}] & \mf{0} \end{pmatrix}^T -\mathbf{U}$ is the source term, and where $M_0[\mathbf{U}]=4\pi M[\mf{U}]$ is defined by (\ref{max}). The conserved variables in the Maxwellian and the in the collision frequency can be related to $\mf{U}$ by $\mf{W}=\int_{\mathbb{R}_+}f_0(v)\boldsymbol\Psi(v)v^2 \mathrm{d}v$. The integration process leads to a last unknown flux $\mathbf{f}_2$. In order to obtain a complete model, this hierarchy needs to be close by expressing this flux in terms of the first two moments $\mathbf{U}$. An indirect way of proceeding is to chose a specific ansatz $\hat{f}[\mf{U}]$ for the distribution function shape and to set:
\begin{equation}
    \forall \boldsymbol\Omega \in \sphere,~f(\boldsymbol\Omega)=\hat{f}[\mathbf{U}](\boldsymbol\Omega),
\end{equation}
where $\hat{f}[\mf{U}]$ is the ansatz which realizes the moments: $\crochet{\mf{m}\hat{f}}=\mf{U}$. According to this definition, the last flux is
\begin{equation}
    \mathbf{f}_2(\mf{U})=\crochet{\boldsymbol\Omega\otimes\boldsymbol\Omega \hat{f}[\mf{U}]}.
\end{equation}
The M1 closure relies on an entropic argument to enforce the distribution shape and to compute this flux as a function of $f_0$ and $\mf{f}_1$. Linked to this closure is the notion of moments realizability:
\begin{defi} [Moment realizability]
A moment vector $\mathbf{U}$ is realizable if there exists a non-negative distribution function $f$ such that $\crochet{\mathbf{m}f}=\mathbf{U}$.
\end{defi}
\begin{prop}
Let $\mathbf{U}=\begin{pmatrix} f_0 & \mf{f}_1\end{pmatrix}^T$ and $u=\frac{||\mf{f}_1||}{f_0}$. The moment vector is realizable if and only if $f_0 > 0$ and $u<1$, or $ \mathbf{U}=\mathbf{0}$.
\end{prop}

\begin{prop}[M1 distribution function] \label{prop:distribm1}
Let $\mathbf{U}=\begin{pmatrix}f_0 & \mathbf{f}_1\end{pmatrix}^T$ be a vector of realizable moments. If $f_0>0$, then the distribution function $\hat{f}[\mf{U}]$ which minimizes the Boltzmann entropy functional $h(f)= \crochet{f\ln{f}-f}$ under the constraint $\crochet{\mathbf{m} \hat{f}}=\mathbf{U}$ is
\begin{equation}
    \forall \Omega \in \sphere,~\hat{f}[\mf{U}](\boldsymbol\Omega)=e^{\mathbf\Lambda\cdot\mathbf{m}(\boldsymbol\Omega)}=f_0 \dfrac{||\boldsymbol\beta||}{4\pi\sinh{||\boldsymbol\beta}||} e^{\boldsymbol\beta \cdot \boldsymbol\Omega},
\end{equation}
where $\mathbf\Lambda = \begin{pmatrix} \alpha & \boldsymbol\beta \end{pmatrix}^T$ is the vector of entropic variables defined through the following relations
\begin{subequations}
    \begin{align}
        &\alpha = \ln{\left(f_0\frac{||\boldsymbol\beta||}{4\pi\sinh ||\boldsymbol\beta||}\right)}, \\
        &||\boldsymbol\beta|| = z^{-1}(u), \\
        &\frac{\boldsymbol\beta}{||\boldsymbol\beta||} =\mf{d},
    \end{align}
\end{subequations}
where $u=\frac{||\mf{f}_1||}{f_0}$ is the anisotropic factor,  $\mf{d}=\frac{\mf{f}_1}{||\mf{f}_1||}$ is the velocity direction and $z(x)= \coth{x}-x^{-1}$ is an invertible odd function in $[0,1]$, continuously extendable at $x=0$.
\end{prop}
\begin{prop}[M1 closure]
The third moment of the M1 distribution function $\hat{f}[\mf{U}]$ is
\begin{equation} \label{q_m1}
    \mf{f}_2=f_0\dfrac{u}{||\boldsymbol\beta||}\mathbf{I}_3 + f_0\left(1-3\dfrac{u}{||\boldsymbol\beta||}\right) \dfrac{\boldsymbol\beta}{||\boldsymbol\beta||} \otimes \dfrac{\boldsymbol\beta}{||\boldsymbol\beta||}.
\end{equation}
\end{prop}
The detailed computation of the moments of the M1 distribution function can be found in \ref{ann1}. System (\ref{model_moment}) closed with relation (\ref{q_m1}) is the M1 model associated with the kinetic equation (\ref{equation_base}). In the case of a zero distribution function the first angular moment $f_0$ is zero and the closing procedure is not applicable anymore because the velocity $u$ and hence $\boldsymbol\beta$ are not well defined anymore. But the continuity of the entropy functional $h$ at $f=0$ allows to set $\hat{f}=0$ and therefore $\mf{f}_2=0$.

This model satisfies several properties. First it can be shown (see \cite{dubroca1999etude} for details) that (\ref{model_moment})-(\ref{q_m1}) is a hyperbolic model and that the moments realizability is ensured. Moreover the macroscopic variables of the M1 model satisfy the same diffusion equations as the kinetic equation. This property shows that this model allows an intermediate description between the full kinetic model and the macroscopic one.


\section{A UGKS based numerical scheme for the M1 model for structured meshes\label{sec:ugks}}
As the M1 model (\ref{model_moment})-(\ref{q_m1}) is a hyperbolic system, a standard approximate Riemann solver could be used to solve it. However, without a special treatment of the source term, this scheme cannot correctly discretize the diffusion equation (\ref{diffusion}) in the corresponding limit. Instead of artificially modifying the solver to recover correct fluxes in this regime, our strategy is to rely on a robust asymptotic preserving scheme for the kinetic equation. The Unified Gas Kinetic Scheme (UGKS) is a good choice since it naturally ensures this property by construction. In this section, we detail the UGKS construction is detailed for the electron transport equation before constructing the new UGKS-based scheme for the M1 moment model.

    \subsection{UGKS}
    \subsubsection{A finite volume formulation}
First, we start with a standard finite volume formulation. Let \newline $C_{i,j}=[x_{i-1/2},x_{i+1/2}]\times[y_{j-1/2},y_{j+1/2}]$ be a control volume of size $\Delta x \Delta y$ with centre $\mathbf{x}_{i,j}$  and $[t_n,t_{n+1}]$ be a time interval of size $\Delta t$. We define the averages of the macroscopic variables and distribution function on cell $(i,j)$ at time $t_n$:
\begin{subequations}
    \begin{align}
        \boldsymbol{W}_{i,j}^n&=\dfrac{1}{\Delta x \Delta y} \iint_{C_{i,j}} \ccrochet{\boldsymbol\Psi f(t_n,x,y,\cdot,\cdot)} \mathrm{d}x\mathrm{d}y, \\
    f_{i,j}^n&=\dfrac{1}{\Delta x \Delta y} \iint_{C_{i,j}} f(t_n,x,y,\cdot,\cdot) \mathrm{d}x\mathrm{d}y.
    \end{align}
\end{subequations}
The finite volume formulations of the kinetic equation is obtained by integrating equations ($\ref{equation_base}$) over the control volume and over the time interval. This formulation emphasizes the evolution of the distribution function average through the cell interface fluxes during a time step:
\begin{equation}\label{finite_volumef}
    \dfrac{f_{i,j}^{n+1}-f_{i,j}^n}{\Delta t} + \dfrac{\phi_{i+1/2,j}-\phi_{i-1/2,j}}{\Delta x} + \dfrac{\phi_{i,j+1/2}-\phi_{i,j-1/2}}{\Delta y}
    =\nu(\boldsymbol{W}_{i,j}^{n+1})\left(M[\boldsymbol{W}_{i,j}^{n+1}]-f_{i,j}^{n+1}\right),
\end{equation}
where the microscopic numerical fluxes across the middle of the eastern and northern interfaces are
\begin{subequations}
    \begin{align}
    \phi_{i+1/2,j}(v,\boldsymbol\Omega) &=  \dfrac{1}{\eta\Delta t} \dfrac{1}{\Delta y}\int_{y_{j-1/2}}^{y_{j+1/2}}\int_{t_n}^{t_{n+1}} v \Omega_x f(t,x_{i+1/2},y,v,\boldsymbol\Omega) \mathrm{d}t \mathrm{d}y, \\
    \phi_{i,j+1/2}(v,\boldsymbol\Omega) &= \dfrac{1}{\eta\Delta t} \dfrac{1}{\Delta x}\int_{x_{i-1/2}}^{x_{i+1/2}}\int_{t_n}^{t_{n+1}} v \Omega_y f(t,x,y_{j+1/2},v,\boldsymbol\Omega) \mathrm{d}t \mathrm{d}x.
    \end{align}
\end{subequations}
For first and second order schemes, we use a midpoint rule to approximate the fluxes along the interfaces by their value on the centre. For example, on the eastern interface we have
\begin{equation}
    \phi_{i+1/2,j}(v,\boldsymbol\Omega) = \dfrac{1}{\eta\Delta t} \int_{t_n}^{t_{n+1}}v\Omega_x f(t,\mf{x}_{i+1/2,j},v,\boldsymbol\Omega)\mathrm{d}t + \mathcal{O}(\Delta y^2).
\end{equation}
Moreover, an implicit approximation of the collision term is chosen to obtain an uniformly stable scheme which leads to an unknown term $\mf{W}_{i,j}^{n+1}$. In order to perform an explicit computation of the distribution function, the macroscopic vector at time $t_{n+1}$ is updated first using a finite volume formulation of the macroscopic conservation law (\ref{conservation_macro}):
\begin{equation} \label{vf_macro}
    \dfrac{\boldsymbol{W}_{i,j}^{n+1}-\boldsymbol{W}_{i,j}^n}{\Delta t} + \dfrac{1}{\Delta x} \left(\boldsymbol\Phi_{i+1/2,j}-\boldsymbol\Phi_{i-1/2,j}\right) + \dfrac{1}{ \Delta y}  \left(\boldsymbol\Phi_{i,j+1/2}-\boldsymbol\Phi_{i,j-1/2}\right) =0,
\end{equation}
where the macroscopic fluxes vectors are moments of the microscopic fluxes:
\begin{subequations}
    \begin{align}
    \boldsymbol\Phi_{i+1/2,j} &= \dfrac{1}{\eta\Delta t} \dfrac{1}{\Delta y}\int_{y_{j-1/2}}^{y_{j+1/2}}\int_{t_n}^{t_{n+1}} \ccrochet{ v\boldsymbol\Psi(v) \Omega_x f(t,x_{i+1/2},y,v,\boldsymbol\Omega} \mathrm{d}t \mathrm{d}y, \\
    \boldsymbol\Phi_{i,j+1/2} &= \dfrac{1}{\eta\Delta t} \dfrac{1}{\Delta x}\int_{x_{i-1/2}}^{x_{i+1/2}}\int_{t_n}^{t_{n+1}} \ccrochet{ v\boldsymbol\Psi(v) \Omega_y f(t,x,y_{j+1/2},v,\boldsymbol\Omega} \mathrm{d}t \mathrm{d}x.
    \end{align}
\end{subequations}
The same midpoint rule is used to approximate the integral over the interface. Developing a finite volume scheme for equation ($\ref{equation_base}$) involves giving a consistent and conservative approximation of the microscopic numerical fluxes and therefore of the macroscopic ones as $\boldsymbol\Phi_{i+1/2,j}=\ccrochet{\boldsymbol\Psi\phi_{i+1/2,j}}$. At this stage, the velocity variable $v$ as well as the direction variable $\boldsymbol\Omega$ are kept continuous and omitted when not needed. In the following, the development of the numerical flux is only detailed on the eastern interface for the sake of simplicity.

        \subsubsection{A characteristic based approach}
The numerical flux of UGKS rests on the integral representation of the kinetic equation solution. That procedure allows to naturally take into account the collision term which induces the diffusion mechanism. Assuming a constant collision frequency $\nu$, equation (\ref{equation_base}) is equivalent to:
\begin{equation}
    \derivee{}{t}\left(e^{\nu t} f(t,\mathbf{x}+\dfrac{v}{\eta}t\boldsymbol{\Omega},v,\boldsymbol\Omega)\right)=\nu e^{\nu t}M[\mathbf{W}](t,\mathbf{x}+\dfrac{v}{\eta}t\boldsymbol\Omega,v,\boldsymbol\Omega)
\end{equation}
However, the collision frequency is \textit{a priori} non-constant as it is a function of the macroscopic state $\mathbf{W}$. We assume that its variations are negligible at the scale of a cell and of a time step to consider this expression as an approximation between two given times $t_n$ and $t>t_n$ around an interface. For example around $\mathbf{x}_{i+1/2,j}$, a
simple integration in time gives
\begin{equation} \label{representation_interface}
    \begin{aligned}
        f(t,\mathbf{x}_{i+1/2,j},v,\boldsymbol\Omega) & \approx e^{-\nu_{i+1/2,j}^n (t-t_n)} f(t_n,\mathbf{x}_{i+1/2,j}-\dfrac{v}{\eta}(t-t_n)\boldsymbol{\Omega},v,\boldsymbol\Omega) \\
        & + \nu_{i,j+1/2}^n \int_{t_n}^t e^{-\nu_{i+1/2,j}^n (t-s)} M[\mathbf{W}](s,\mathbf{x}_{i+1/2,j}-\dfrac{v}{\eta}(t-s)\boldsymbol{\Omega},v) \mathrm{d}s,\\
    \end{aligned}
\end{equation}
where the mean collision frequency is $\nu_{i+1/2,j}=\frac{1}{2}(\nu_{i,j}+\nu_{i+1,j})$. 
In order to evaluate the horizontal numerical flux $\phi_{i+1/2,j}$ from relation ($\ref{representation_interface}$), cell average based reconstructions in space and time of $f$ and $M[f]$ need to be introduced. Appropriate choices are mandatory to preserve the asymptotics and achieve second order convergence in space. The Maxwellian reconstruction $\tilde M$ is chosen to be continuous at the interface and defined as
 \begin{equation} \label{reconstruction_maxwellian}
    \tilde{M}(t,x,y,\cdot)= 
    M[\mf{W}_{i+1/2,j}^n]+
    \begin{cases}
    {\delta^L_x} M_{i+1/2,j}^n ({x}-{x}_{i+1/2})& \text{ if } {x} \leq x_{i+1/2} \\
    {\delta^R_x} M_{i+1/2,j}^n ({x}-{x}_{i+1/2})& \text{ if } {x} \geq x_{i+1/2} \\
    \end{cases},
\end{equation}
where $M[\mf{W}_{i+1/2,j}^n]$ is the Maxwellian at the interface macroscopic state, given by the arithmetic mean $\mf{W}_{i+1/2,j}^n=\frac{1}{2}(\mf{W}_{i,j}^n+\mf{W}_{i+1,j}^n)$ and $\delta_x^{LR} M_{i+1/2,j}^n$ are the Maxwellian spatial slopes. First, it should be noted that $\tilde{M}$ is no longer a Maxwellian. Secondly, the tangential part of the reconstruction is neglected as it does not contribute to the macroscopic fluxes and as a consequence does not alter the asymptotic behaviour. The spatial derivative of the Maxwellian (given by the chain rule) allows to compute the slopes in terms of the derivatives of the conservative variables:
\begin{equation} \label{dl_max}
\begin{aligned}
\partial_{x} M[\mathbf{W}]&= \left(\partial_{x} \rho\right)\partial_\rho M[\mathbf{W}](v)+\left(\partial_{x} q\right) \partial_q M[\mathbf{W}](v),  \\
&= \boldsymbol\Psi(v)^T\mf{K}[\mf{W}] \partial_x \mf{W} , \\
\end{aligned}
\end{equation}
where the matrix $\mf{K}$ is
\begin{equation} \label{mat_k}
    K(\mf{W})=M[\mf{W}]\begin{pmatrix}
    \dfrac{5}{2\rho} & -\dfrac{3}{2q} \\
    -\dfrac{3}{2q} & \dfrac{3\rho}{2q^2} \\
    \end{pmatrix}.
\end{equation}
Finally, we approximate $\partial_x \mf{W}$ in (\ref{dl_max}) by finite difference formula to define the left and right slopes $\partial_x^{LR}M_{i+1/2,j}^n$ of the Maxwellian reconstruction (\ref{reconstruction_maxwellian}):
\begin{subequations} \label{pentes_maxwellienne}
\begin{align}
{\delta^{L}_x} M_{i+1/2,j}^n&=\boldsymbol\Psi^T \mf{K}(\mf{W}_{i+1/2,j}^n) \dfrac{\mf{W}_{i+1/2,j}^n-\mf{W}_{i,j}^n}{\Delta x/2}, \\
{\delta^{R}_x} M_{i+1/2,j}^n&=\boldsymbol\Psi^T \mf{K}(\mf{W}_{i+1/2,j}^n) \dfrac{\mf{W}_{i+1,j}^n-\mf{W}_{i+1/2,j}^n}{\Delta x/2}.
\end{align}
\end{subequations}
The constructed distribution function $\tilde{f}$ is defined by a more conventional piecewise continuous affine function:
\begin{equation} \label{reconstruction_f}
   \tilde{f}(t_n,x,y,\cdot,\cdot)=\begin{cases}
    f_{i,j}^n + {\delta_x}f_{i,j}^n ({x}-{x}_{i})& \text{ if } {x}\leq x_{i+1/2} \\
    f_{i+1,j}^n + {\delta_x}f_{i+1,j}^n ({x}-{x}_{i+1})& \text{ if } {x}\geq x_{i+1/2} \\
    \end{cases},
\end{equation}
where $\delta_x f_{i,j}^n$ is the distribution function slope. The tangential part of the reconstruction is also neglected for the sake of simplicity and for computational reasons as this choice leads to an enlarged stencil. The slope needs to be limited to ensure the decrease of the total variation. Let $\psi$ be a TVD slope limiter (for example the van Leer limiter is $\psi(x,y)=(\sign(x)+\sign(y)) \frac{|x||y|}{|x|+|y|}$ ). Then the limited slope is given by a simple finite difference formula:
\begin{equation}
    \delta_xf_{i,j}^n=\psi\left(\dfrac{f_{i+1,j}^n-f_{i,j}^n}{\Delta x},\dfrac{f_{i,j}^n-f_{i-1,j}^n}{\Delta x}\right).
\end{equation}
To evaluate $\phi_{i+1/2,j}$, the reconstructed quantities $\tilde{f}$ and $\tilde{M}$ are employed in ($\ref{representation_interface}$) before time integration. It should be noted that in the diffusion limit, the foot of the characteristics might be arbitrarily far from the interface. However, due to the collision mechanism, the particles are constrained near the interface (as shown by the exponential term in ($\ref{representation_interface}$)). Therefore, it is legitimate to neglect the influence of remote particles by extending the reconstructions validity domain. Finally, the microscopic numerical flux takes the following form
\begin{equation} \label{flux_micro}
\begin{split}
\phi_{i+1/2,j}(v,\boldsymbol\Omega)=&A_{i+1/2,j}^n v \Omega_x \left(\indi{f_{i,j}^{n(+)}}{f_{i+1,j}^{n(-)}}{\Omega_x}\right)\\
&+ B_{i+1/2,j}^n v^2 \Omega_x^2(\delta_x f_{i,j}^n \mathbb{1}_{\Omega_x>0}+\delta_x f_{i+1,j}^n \mathbb{1}_{\Omega_x<0})\\
 &+  C_{i+1/2,j}^n v \Omega_x M[\mf{W}_{i+1/2,j}^n]\\
 &+  D_{i+1/2,j}^n v^2 \Omega_x^2 \left(\indi{{\delta^L_x }M_{i+1/2,j}^n}{{\delta^R_x }M_{i+1/2,j}^n}{\Omega_x}\right),
\end{split}
\end{equation}
and the macroscopic ones are
\begin{subequations} \label{ugks_macro}
    \begin{align}
    \begin{split}
    \Phi^\rho_{i+1/2,j}=&A_{i+1/2,j}^n \ccrochet{\indi{v\Omega_x f_{i,j}^{n(+)}}{v\Omega_xf_{i+1,n}^{n(-)}}{\Omega_x} } \\
    &+ B_{i+1/2}^n \ccrochet{v^2 \Omega_x^2 \delta_x f_{i,j}^n \mathbb{1} _{\Omega_x>0} + v^2\Omega_x^2\delta_x f_{i+1,j}^n \mathbb{1}_{\Omega_x<0}}\\
    &+\dfrac{2D_{i+1/2,j}^n}{3\Delta x} (q_{i+1,j}^n-q_{i,j}^n),
    \end{split} \\
    \begin{split}
    \Phi^q_{i+1/2,j}=&\dfrac{A_{i+1/2,j}^n}{2} \ccrochet{\indi{v^3\Omega_x f_{i,j}^{n(+)}}{v^3\Omega_x f_{i+1,n}^{n(-)}}{\Omega_x} } \\
    &+ \dfrac{B_{i+1/2,j}^n}{2} \ccrochet{v^4 \Omega_x^2 \delta_x f_{i,j}^n \mathbb{1} _{\Omega_x>0} + v^4\Omega_x^2\delta_x f_{i+1,j}^n \mathbb{1}_{\Omega_x<0}}\\
    &+\dfrac{2D_{i+1/2,j}^n}{3\Delta x} \dfrac{5(q_{i+1/2,j}^n)^2}{3\rho_{i+1/2,j}^n}\left( -
    \dfrac{\rho_{i+1,j}^n-\rho_{i,j}^n}{\rho_{i+1/2,j}^n} + 2 \dfrac{q_{i+1,j}^n-q_{i,j}^n}{q_{i+1/2,j}^n}\right).
    \end{split}
    \end{align}
\end{subequations}
where $f_{i,j}^{n(\pm)}=f_{i,j}^n\pm\frac{\Delta x}{2}\delta_xf_{i,j}^n$. The detailed calculation of the macroscopic fluxes can be found in \ref{ann0}.
The integration coefficients $A_{i+ 1/2,j}$, $B_{i+ 1/2,j}$, $C_{i+ 1/2,j}$ and $D_{i+ 1/2,j}$ are interface values of functions
\begin{subequations} \label{coeff}
\begin{align}
A(\Delta t,\eta,\epsilon,\sigma)&=\dfrac{-1}{\eta}\dfrac{(1-e^w)}{w}, \\
B(\Delta t,\eta,\epsilon,\sigma)&=\dfrac{1}{\sigma}\dfrac{\epsilon}{\eta}  \left( e^{w} +\dfrac{1-e^w}{w} \right), \\
C(\Delta t,\eta,\epsilon,\sigma)&=\dfrac{1}{\eta}\left(1+\dfrac{1-e^w}{w}\right),\\    
D(\Delta t,\eta,\epsilon,\sigma)&=\dfrac{-1}{\sigma}\dfrac{\epsilon}{\eta} \left(1+e^{w}+2\dfrac{1-e^w}{w}\right),
\end{align}
\end{subequations}
at $\sigma_{i+1/2,j}=\frac{\sigma_{i,j}+\sigma_{i+1,j}}{2}$ and where $w=-\nu \Delta t=-\frac{\sigma}{\epsilon\eta}\Delta t$.
        
    \subsubsection{Asymptotic behaviour}
The behavior of the scheme is examined in both the free transport regime ($\epsilon \to +\infty$) and in the diffusion regime ($\epsilon \to 0$) at the diffusion scale $\eta=\epsilon$. In this asymptotic analysis, the opacity $\sigma(\mathbf{W}(\mathbf{x}))$ is assumed bounded. By studying the limits of the integration coefficients, one can show that the macroscopic fluxes in the diffusion limit are
\begin{subequations} \label{flux_diffusion}
    \begin{align}
    \begin{split}
    \Phi_{i+1/2,j}^\rho&\underset{\epsilon \to 0}{\longrightarrow}
       \dfrac{-2}{3\sigma_{i+1/2,j}} \dfrac{q_{i+1,j}^n-q_{i,j}^n}{\Delta x},
    \end{split} \\
    \begin{split}
    \Phi^q_{i+1/2,j}&\underset{\epsilon \to 0}{\longrightarrow}\dfrac{-2}{3\sigma_{i+1/2,j}} \dfrac{5(q_{i+1/2,j}^n)^2}{3\rho_{i+1/2,j}^n}\left( \dfrac{-1}{\rho_{i+1/2,j}^n}
    \dfrac{\rho_{i+1,j}^n-\rho_{i,j}^n}{\Delta x} + \dfrac{2}{q_{i+1/2,j}^n} \dfrac{q_{i+1,j}^n-q_{i,j}^n}{\Delta x}\right).
    \end{split}
    \end{align}
\end{subequations}
These fluxes correctly discretize diffusion equation (\ref{diffusion}) in the limit regime. A similar analysis can be performed in the free transport regime at the fixed transport scale $\eta$ and with $\epsilon \to +\infty$. The microscopic flux is
\begin{equation}
    \label{flux_micro_limit}
    \phi_{i+1/2,j} \underset{\epsilon \to \infty}{\longrightarrow} \dfrac{v\Omega_x}{\eta} (f_{i,j}^{n(+)} \mathbb{1}_{\Omega_x>0} + f_{i+1,j}^{n(-)} \mathbb{1}_{\Omega_x<0}) 
    -\Delta t \dfrac{v^2\Omega_x^2}{2\eta^2} (\delta_x f_{i,j}^n \mathbb{1}_{v>0} + \delta_x f_{i+1,j}^n \mathbb{1}_{v<0}),
\end{equation}
which is the standard second order scheme for the linear transport equation \cite{leveque2002finite}.

    \subsection{UGKS-M1}
    \subsubsection{A M1 closure of the UGKS}

The new scheme for the M1 moment model of electron transport is derived using the same approach as explained in our previous work \cite{feugeas2024asymptotic}. The main idea is to apply the UGKS to the M1 distribution function $(\hat{f}_{i,j}^n(v,\boldsymbol\Omega))_{i,j}$ reconstructed from the moments $(\mathbf{U}_{i,j}^n(v))_{i,j}$. Then, the macroscopic variables at time $t_{n+1}$ are the moments of $(f_{i,j}^{n+1}(v,\boldsymbol\Omega))_{i,j}$.
In other words, a finite volume formulation of the M1 model (\ref{model_moment})-(\ref{q_m1}) is obtained by integrating the kinetic equation formulation (\ref{finite_volumef}) against $\mf{m}(\boldsymbol\Omega)$:

\begin{equation}\label{ugsk_m1_fv}
    \dfrac{\mathbf{U}_{i,j}^{n+1}-\mathbf{U}_{i,j}^{n}}{\Delta t} + \dfrac{(\boldsymbol{\chi}_{i+1/2,j}-\boldsymbol{\chi}_{i-1/2,j})}{\Delta x}+ \dfrac{(\boldsymbol{\chi}_{i,j+1/2}-\boldsymbol{\chi}_{i,j-1/2})}{\Delta y}=\nu_{i,j}^{n+1}\mathbf{S}(\mathbf{U}_{i,j}^{n+1}),
\end{equation}
where $\mf{U}_{i,j}^{n}=\begin{pmatrix}f_{i,j}^{0,n} & f_{i,j}^{1,n}\end{pmatrix}^T$ are the M1 variables, \\ $\boldsymbol\chi_{i+1/2,j}=\crochet{\mf{m}(\boldsymbol\Omega)\phi_{i+1/2,j}(\cdot,\boldsymbol\Omega)}=\begin{pmatrix}\chi^0_{i+1/2,j} & \boldsymbol\chi_{i+1/2,j}^1\end{pmatrix}^T \in \mathbb{R}^4$ is the M1 numerical flux at the mesoscopic scale and $\mf{S}(\mf{U}_{i,j}^{n+1})=\begin{pmatrix}M_0[\mf{U}_{i,j}^{n+1}] & 0\end{pmatrix}^T$ is the source term. This flux vector is computed by integrating the microscopic UGKS flux $\phi_{i+1/2,j}$ (\ref{flux_micro}) against $\mf{m}(\boldsymbol\Omega)$ with the closed distribution function: $f_{i,j}^n=\hat{f}[\mf{U}_{i,j}^n]$, that will simply be denoted by $\hat{f}_{i,j}^n$ in the following. First, we define the fluxes without the first order term in (\ref{reconstruction_f}) (the distribution function reconstruction is constant per cell). Thus, the mesoscopic fluxes are
\begin{subequations} \label{flux_micro1_ugksm1}
    \begin{align}
    \begin{split}
    \chi^0_{i+1/2,j}(v)=&A_{i+1/2,j}^n v \crochet{\Omega_x \hat{f}_{i,j}^{n} \mathbb{1}_{\Omega_x> 0}+\Omega_x\hat{f}_{i+1,j}^{n} \mathbb{1}_{\Omega_x< 0}} \\
    &+\dfrac{D_{i+1/2,j}^n}{3\Delta x} v^2 \boldsymbol\Psi^T \mf{K}_0\left(\mf{W}_{i+1/2,j}^n\right) \left(\mf{W}_{i+1,j}^n-\mf{W}_{i,j}^n \right),
    \end{split} \\
    \begin{split}
    \boldsymbol\chi^1_{i+1/2,j}(v)=&A_{i+1/2,j}^n v \crochet{\Omega_x\boldsymbol\Omega \hat{f}_{i,j}^{n} \mathbb{1}_{\Omega_x> 0}+\Omega_x\boldsymbol\Omega\hat{f}_{i+1,j}^{n} \mathbb{1}_{\Omega_x< 0}} \\
    &+\dfrac{C_{i+1/2,j}^n}{3} v M_0[\mf{W}_{i+1/2,j}^n] \mathbf{e}_x \\
    &-\dfrac{D_{i+1/2,j}^n}{4\Delta x} v^2\boldsymbol\Psi^T \mf{K}_0\left(\mf{W}_{i+1/2,j}^n\right) \left(\mf{W}_{i+1,j}^n-2\mf{W}_{i+1/2,j}^n+\mf{W}_{i,j}^n \right)\mf{e}_x,
    \end{split}
    \end{align}
\end{subequations}
where $\mf{e}_x$ denotes the unit vector in the x-direction. As in UGKS, the source term in the finite volume formulation (\ref{ugsk_m1_fv}) is implicit and requires updating the macroscopic variables first. To achieve this, the same macroscopic finite volume formulation (\ref{vf_macro}) is used. In that case, the first order macroscopic fluxes are the UGKS ones (\ref{ugks_macro}) where the linear part of the distribution function reconstruction is omitted and where $f_{i,j}^n=\hat{f}_{i,j}^n$:
\begin{subequations}\label{flux_macro1_ugksm1}
    \begin{align} 
    \begin{split}
    \Phi_{i+1/2,j}^\rho=&A_{i+1/2,j}^{n} \ccrochet{\indi{v\Omega_x \hat{f}_{i,j}^{n}}{v\Omega_x\hat{f}_{i+1,n}^{n}}{\Omega_x} } \\
    &+\dfrac{2D_{i+1/2,j}^n}{3\Delta x} (q_{i+1,j}^n-q_{i,j}^n), 
    \end{split} \\
    \begin{split}
    \Phi_{i+1/2,j}^q=&\dfrac{A_{i+1/2,j}^n}{2} \ccrochet{\indi{v^3\Omega_x \hat{f}_{i,j}^{n}}{v^3\Omega_x \hat{f}_{i+1,n}^{n}}{\Omega_x} } \\
    &+\dfrac{2D_{i+1/2,j}^n}{3\Delta x} \dfrac{5(q_{i+1/2,j}^n)^2}{3\rho_{i+1/2,j}^n}\left( -
    \dfrac{\rho_{i+1,j}^n-\rho_{i,j}^n}{\rho_{i+1/2,j}^n} + 2 \dfrac{q_{i+1,j}^n-q_{i,j}^n}{q_{i+1/2,j}^n}\right).
    \end{split}
    \end{align}
\end{subequations}
To sum up, UGKS-M1 consists of the macroscopic part (\ref{vf_macro})-(\ref{flux_macro1_ugksm1}) and the mesoscopic part (\ref{ugsk_m1_fv})-(\ref{flux_micro1_ugksm1}). In order to make this scheme usable and implementable, it is necessary to introduce a speed discretization and a quadrature formula to compute the macroscopic fluxes. Moreover, expressing the half sphere moments of the M1 distribution functions is also required. Finally, restoring second order accuracy can be achieved by using the full distribution function reconstruction in UGKS. These various aspects of the scheme are discussed in the following sections.

\subsubsection{Speed discretization}
Until now, the speed variable $v$ has been kept continuous in the development of UGKS and UGKS-M1. In practice, however, this variable is discretized. One issue is that this variable takes values in an infinite space that needs to be truncated. Depending on the specific physical situation, it is possible to define an upper bound $v_m$ on the speed above which the equilibrium distributions are all below a certain threshold. This maximum speed is a function of the maximum temperature in the domain. Let $(v_k)_{k\leq m}$ be the discrete speeds in $[0,v_m]$ and $\Delta v$ be a speed step. The microscopic variables are updated at each discretization point. On the other hand, macroscopic fluxes require evaluating integrals in speed space. Let $\hat{f}_{i,j,k}^n(\boldsymbol\Omega)$ be the discrete M1 distribution function at speed $v_k$. These integrals are approximated using the composite trapezoidal quadrature formula:
    \begin{equation}
        \forall p\geq 1,~\ccrochet{v^p\hat{f}_{i,j}^n(v,\cdot)} \simeq \dfrac{\Delta v}{2} \sum_{k=0}^{m-1} \left( \crochet{v^p_k \hat{f}_{i,j,k}^n}+\crochet{v^p_{k+1} \hat{f}_{i,j,k+1}^n}\right),
    \end{equation}
where $\Delta v$ is the speed step. The quadrature rule can be chosen in order to exactly integrate the Maxwellian equilibrium. However, this property is not required to recover a good diffusion scheme. 
    
    \subsubsection{Dirichlet boundary condition}
In this section, we describe how Dirichlet boundary conditions can be treated in UGKS and therefore in UGKS-M1. For example in the Western border, the incoming particles distribution is imposed to a Maxwellian  $M[\mf{W}_W(t,\mathbf{x})](v)$ at a given macro state. In order to obtain the numerical flux at this boundary, the distribution function representation is modified to take into account the incoming particles from the boundary:
    \begin{equation}
    \label{representation_interface_bord}
    f(t,\mathbf{x}_{1/2,j},v,\boldsymbol\Omega)=
    \begin{cases}
    M[\mf{W}_W(t,\mathbf{x})](v) & \text{ if } \Omega_x>0 \\
    \begin{array}{l}
e^{-\nu^n_{1/2,j}(t-t_n)} f(t_n,\mathbf{x}_{1/2,j}-\dfrac{v}{\eta}(t-t_n)\boldsymbol\Omega,v,\boldsymbol\Omega) \\
+\nu^n_{1/2,j}\displaystyle \int_{t_n}^t  e^{-\nu^n_{1/2,j}(t-s)} M[f](s,\mf{x}_{1/2,j} - \dfrac{v}{\eta}(t-s)\boldsymbol\Omega,v) \mathrm{d}s 
\end{array}& \text{ if } \Omega_x<0
    \end{cases}.
\end{equation}
Using similar distribution function and Maxwellian reconstructions, the first order microscopic flux is
\begin{equation} \label{flux_micro_boundary}
\begin{aligned}
\phi_{1/2,j}(v,\boldsymbol\Omega)&=\dfrac{v\Omega_x}{\eta}M[\mf{W}_W(t_n,\mathbf{x}_{1/2,j})](v) \mathbb{1}_{\Omega_x > 0} \\
&+\left(A_{1/2,j}^n v \Omega_x {f_{i+1,j}^{n}}+  C_{1/2,j}^n v \Omega_x M[\mf{W}_{1/2,j}^n] + D_{1/2,j}^n v^2 \Omega_x^2 {{\delta^R_x }M_{i+1/2,j}^n} \right) \mathbb{1}_{\Omega_x < 0}.
\end{aligned}
\end{equation}
In order to ensure a correct asymptotic behaviour of the macroscopic boundary flux in the diffusion limit, the interface macroscopic variable vector $\mf{W}^n_{1/2,j}$ is set to
\begin{equation}
    \mf{W}_{1/2,j}^n=\mf{W}_W(t_n,\mathbf{x}_{1/2,j}).
\end{equation}
Finally, the macroscopic and the UGKS-M1 boundary fluxes are moments of this flux.
    
    \subsubsection{Half sphere moments evaluation \label{sec:quad}}
Despite the fact that the moments of the M1 distribution function on the unit sphere are analytical, the half moments are not. In this section we propose a numerical procedure to evaluate the required integrals in (\ref{flux_micro1_ugksm1}) for UGKS-M1. To begin with, the half moments are given over the half sphere $\{ \boldsymbol\Omega \in \sphere ~|~\Omega_z \geq 0 \}$. First, the first two half moments are reformulated in the form of 1D integrals:
\begin{equation}
    \begin{aligned}
    \forall k\geq0,~\crochet{\Omega_z^k \hat{f}\mathbb{1}_{\Omega_z\geq 0}}&=\int_{\sphere} \Omega_z^ke^{\alpha+\boldsymbol\beta \cdot \boldsymbol\Omega}\mathbb{1}_{\Omega_z\geq 0}\mathrm{d}\Omega, \\
    &=f_0\dfrac{||\boldsymbol\beta||}{4\pi \sinh{||\boldsymbol\beta||}}\int_{0}^1\mu^k e^{\beta_z \mu}\int_{-\pi}^\pi e^{\beta_{xy}\sqrt{1-\mu^2}\cos{(\theta-\theta_0)}}\mathrm{d}\phi\mathrm{d}\mu, \\
    &=f_0\dfrac{||\boldsymbol\beta||}{2 \sinh{||\boldsymbol\beta||}}\int_{0}^1\mu^k e^{\beta_z \mu} \mathrm{I}_0(\beta_{xy}\sqrt{1-\mu^2})\mathrm{d}\mu,
    \end{aligned}
\end{equation}
where $\beta_{xy}=\sqrt{\beta_x^2+\beta_y^2}$ is the norm of the projection of vector $\boldsymbol\beta$ on the plane $\{z=0\}$, $\phi_0=\arctan{\frac{\beta_y}{\beta_x}}-(1-\sign{\beta_x})\frac{\pi}{2}$ is the phase and $\mathrm{I}_n$ is the modified Bessel function of the first kind defined as
\begin{equation}
    \forall n\in \mathbb{N},~\forall x\in \mathbb{R}_+,~\mathrm{I}_n(x)=\dfrac{1}{\pi}\int_0^\pi \cos{(n\phi)}e^{x\cos{\phi}}\mathrm{d}\phi.
\end{equation}
The other moments can be expressed similarly:
\begin{equation}
    \crochet{\Omega_z\begin{pmatrix} \Omega_x \\ \Omega_y \end{pmatrix} \hat{f}\mathbb{1}_{\Omega_z>0}}=f_0\dfrac{||\boldsymbol\beta||}{2 \sinh{||\boldsymbol\beta||}}\begin{pmatrix} \cos{\phi_0} \\ \sin{\phi_0} \end{pmatrix}\int_{0}^1\mu\sqrt{1-\mu^2} e^{\beta_z \mu} \mathrm{I}_1(\beta_{xy}\sqrt{1-\mu^2})\mathrm{d}\mu.
\end{equation}
Then, a Gauss-Legendre quadrature is performed on all of those 1D integrals and the approximated integrals are denoted by $\crochet{\cdot}_{GL}$. This procedure leads to inconsistent fluxes in the transport regime. To recover this property, the half moments formulas needs to be consistent with the total value on the unit sphere. For example, the first moment should verify:
\begin{equation}
    \crochet{\Omega_z \hat{f}\mathbb{1}_{\Omega_z> 0}}_{GL}+\crochet{\Omega_z \hat{f}\mathbb{1}_{\Omega_z< 0}}_{GL}=\crochet{\Omega_z \hat{f}}.
\end{equation}
As the total value of this integral is known, the half moments approximations can be renormalised (when needed):
\begin{equation}
    \crochet{\Omega_z \hat{f}\mathbb{1}_{\Omega_z> 0}} \simeq 
    \dfrac{\crochet{\Omega_z\hat{f}}}{\crochet{\Omega_z \hat{f}\mathbb{1}_{\Omega_z\geq 0}}_{GL}+\crochet{\Omega_z \hat{f}\mathbb{1}_{\Omega_z\leq 0}}_{GL}} \crochet{\Omega_z \hat{f}\mathbb{1}_{\Omega_z> 0}}_{GL},
\end{equation}
to recover consistent fluxes. The procedure is similar for the other moments. \\

In order to recover all the required quantities to compute the fluxes, a change of variable is performed. Let $\mf{n}$ be a unit vector of any direction and $\mf{R}$ be the rotation matrix such that $\mf{R}\mf{n}=\mf{e}_z$. This matrix is the composition of two rotation matrix on two different axis:
\begin{equation}
    \mf{R}=
    \begin{pmatrix} \cos\theta & 0 & -\sin\theta \\ 0 & 1 & 0 \\ \sin\theta & 0 & \cos\theta \end{pmatrix}
    \begin{pmatrix} \cos\phi & \sin\phi & 0 \\ -\sin\phi & \cos\phi & 0 \\ 0 & 0 & 1 \end{pmatrix},
\end{equation}
where $(\theta,\phi)$ are the components of $\mf{n}$ in the spherical coordinate system. For the first order scheme, all the required half-moments can be written in the following form:
\begin{equation}
\begin{aligned}
    \crochet{\Omega_n\Omega_i \hat{f}\mathbb{1}_{\Omega_n}}    &=\int_{\mathscr{S}^2} \Omega_z^{'} (\boldsymbol\Omega^{'}\cdot \mf{R}\mf{e}_i) e^{\alpha +\beta^{'}\cdot\boldsymbol\Omega^{'}} \mathbb{1}_{\Omega_z^{'}>0}\mathrm{d}\Omega^{'}.
\end{aligned}
\end{equation}
where $i\in\{x,y,z\}$, $\boldsymbol\Omega^{'}=\mf{R}\boldsymbol\Omega$ and $\boldsymbol\beta^{'}=\mf{R}\boldsymbol\beta$. This integral can be separated into three different integrals: 
\begin{equation}
    \begin{aligned}
            \crochet{\Omega_n\Omega_i \hat{f}\mathbb{1}_{\Omega_n>0}}
            &= \left(\mf{R}\mf{e}_i \cdot \mf{e}_x\right) \int_{\mathscr{S}^2} \Omega_z\Omega_x e^{\alpha + \boldsymbol\beta^{'}\cdot\boldsymbol\Omega} \mathbb{1}_{\Omega_z>0}\mathrm{d}\Omega \\
            &+ \left(\mf{R}\mf{e}_i \cdot \mf{e}_y\right) \int_{\mathscr{S}^2} \Omega_z\Omega_y e^{\alpha + \boldsymbol\beta^{'}\cdot\boldsymbol\Omega} \mathbb{1}_{\Omega_z>0}\mathrm{d}\Omega \\
            &+ \left(\mf{R}\mf{e}_i \cdot \mf{e}_z\right) \int_{\mathscr{S}^2} \Omega_z\Omega_z e^{\alpha + \boldsymbol\beta^{'}\cdot\boldsymbol\Omega} \mathbb{1}_{\Omega_z>0}\mathrm{d}\Omega,
    \end{aligned}
 \end{equation}
 each of them can be evaluated using the method presented above.

        \subsubsection{Second order extensions}
The first order in space UGKS-M1 scheme was obtained by dropping the linear part of the distribution function reconstruction. This term is problematic as the non-linearity introduced by the slope limiter prevents computing analytically the half-moment of the distribution function slopes. In order to achieve a second order convergence rate in space, a different reconstruction of the distribution function is used as proposed in \cite{kxu2010ugks} for the Boltzmann equation of rarefied gas dynamics. First, the vector of conservative variables is reconstructed in both directions. In the horizontal one, the reconstruction is
\begin{equation}
    \mathbf{U}_{i,j}^n (x)=\begin{cases}
    \mathbf{U}_{i,j}^n + \boldsymbol\delta_x \mathbf{U}_{i,j}^n(x-x_{i}) &\text{ if } x<x_{i+1/2}\\
    \mathbf{U}_{i+1,j}^n + \boldsymbol\delta_x \mathbf{U}_{i+1,j}^n(x-x_{i+1}) &\text{ if } x>x_{i+1/2}
    \end{cases},
\end{equation}
where the finite difference slope is $\boldsymbol\delta_x \mathbf{U}_{i,j}^n=\frac{1}{\Delta x}(\mathbf{U}_{i+1,j}^n-\mathbf{U}_{i,j}^n) \mathbf{\phi}(\mathbf{r}_{i,j})$, $\phi$ is a slope limiter and $\mathbf{r}_{i,j}=\frac{\mathbf{U}_{i,j}-\mathbf{U}_{i-1,j}}{\mathbf{U}_{i+1,j}-\mathbf{U}_{i,j}}$ is the local slope defined component-wise. Then, we expand $\hat{f}(\mathbf{U}_{i,j}^n(x))=\exp{(\mathbf{\Lambda}(\mathbf{U}_{i,j}^n(x))\cdot \mathbf{m})}$ in Taylor series (for example when $x<x_{i+1/2}$):
\begin{equation}
    \begin{aligned}
    \hat{f}(\mathbf{U}_{i,j}^n (x))&=\hat{f}(\mathbf{U}_{i,j}^n) + \dfrac{\mathbf{\mathrm{d}\hat{f}}}{\mathrm{d}\mathbf{U}}(\mathbf{U}_{i,j}^n) \cdot \boldsymbol\delta_x \mathbf{U}_{i,j}^n (x-x_{i}) \\
    &=\hat{f}(\mathbf{U}_{i,j}^n) + \mathbf{J}_\mathbf\Lambda (\mathbf{U}_{i,j}^n)^T \mathbf{m}\hat{f}(\mathbf{U}_{i,j}^n)    \cdot \boldsymbol\delta_x \mathbf{U}_{i,j}^n (x-x_{i}),
    \end{aligned}
\end{equation}
where the Jacobian matrix of the transformation is (see \ref{ann4} for details)
\begin{equation} \label{jacob}
\begin{aligned}
    \mathbf{J}_\mathbf\Lambda (\mathbf{U})=\dfrac{f_0^{-1}}{\xi}
     \begin{pmatrix} 1-2\dfrac{u}{||\boldsymbol\beta||} & - \mathbf{u}^T \\ 
                                -\mathbf{u} & \dfrac{||\boldsymbol\beta||}{u}\xi \mathbf{I}_3+(1-\dfrac{||\boldsymbol\beta||}{u}\xi)\dfrac{\mathbf{u}}{||\mathbf{u}||}\otimes\dfrac{\mathbf{u}}{||\mathbf{u}||}\\
    \end{pmatrix}
    \end{aligned}
\end{equation}
and where $\xi=1-2\dfrac{u}{||\boldsymbol\beta||}-u^2$. Finally, the M1 distribution function reconstruction is defined as
\begin{equation}
   f(t_n,\mathbf{x},\cdot,\cdot)=\begin{cases}
    \hat{f}_{i,j}^n + {\delta_x}\hat{f}_{i,j}^n ({x}-{x}_{i,j})& \text{ si } \mathbf{x}\cdot \mathbf{e}_x\leq x_{i+1/2} \\
    \hat{f}_{i+1,j}^n + {\delta_x}\hat{f}_{i+1,j}^n ({x}-{x}_{i+1,j})& \text{ si } \mathbf{x}\cdot \mathbf{e}_x\geq x_{i+1/2} \\
    \end{cases},
\end{equation}
where the slope is 
\begin{equation}
{\delta_x}\hat{f}_{i,j}^n= \mathbf{J}_\mathbf\Lambda(\mathbf{U}_{i,j}^n) \mathbf{\delta U}_{i,j}^n \cdot \mathbf{m}\hat{f}(\mathbf{U}_{i,j}^n).
\end{equation}
The second order fluxes are obtained by integrating the second order UGKS fluxes where the distribution function slopes are computed with the above formula. This procedure leads to two new terms in the microscopic fluxes (\ref{flux_micro1_ugksm1}):
\begin{equation}
    \begin{split}
        ~ & A_{i+1/2,j}^n v \dfrac{\Delta x}{2} \crochet{\Omega_x \mf{m}(\boldsymbol\Omega) {\delta_x}\hat{f}_{i,j}^n \mathbb{1} _{\Omega_x>0} - \Omega_x \mf{m}(\boldsymbol\Omega){\delta_x}\hat{f}_{i+1,j}^n \mathbb{1}_{\Omega_x<0}} \\
        &+ B_{i+1/2}^n v^2\crochet{ \Omega_x^2\mf{m}(\boldsymbol\Omega) {\delta_x}\hat{f}_{i,j}^n \mathbb{1} _{\Omega_x>0} + \Omega_x^2\mf{m}(\boldsymbol\Omega){\delta_x}\hat{f}_{i+1,j}^n \mathbb{1}_{\Omega_x<0}}, \\
    \end{split}
\end{equation}
as well as in the macroscopic ones (\ref{flux_macro1_ugksm1}):
\begin{equation}
    \begin{split}
        ~ & A_{i+1/2,j}^n \dfrac{\Delta x}{2} \ccrochet{v\boldsymbol\Psi(v)\Omega_x  {\delta_x}\hat{f}_{i,j}^n \mathbb{1} _{\Omega_x>0} - v\boldsymbol\Psi(v)\Omega_x {\delta_x}\hat{f}_{i+1,j}^n \mathbb{1}_{\Omega_x<0}} \\
        &+ B_{i+1/2}^n \ccrochet{ v^2\boldsymbol\Psi(v)\Omega_x^2 {\delta_x}\hat{f}_{i,j}^n \mathbb{1} _{\Omega_x>0} + v^2\boldsymbol\Psi(v)\Omega_x^2{\delta_x}\hat{f}_{i+1,j}^n \mathbb{1}_{\Omega_x<0}}. \\
    \end{split}
\end{equation}
The expressions of the higher order half-moments can be found in \ref{ann2}. A re-normalisation step can also be performed on those half-moments to ensure the consistency of the numerical scheme. Because this modification does not seem to alter the solution and significantly increases the computational cost since it requires evaluating $\mf{f}_3$ and $\mf{f}_4$, it is not used in practice. Finally, the full expression of the second order fluxes are written in \ref{ann5}.

To achieve a second order convergence rate in time, the Crank-Nicolson method can be used. Thus, a trapezoidal rule is performed on the source term, which yields to the subsequent finite volume formulation:
\begin{equation}\label{ugsk_m1_fv_order2}
\begin{aligned}
    \dfrac{\mathbf{U}_{i,j}^{n+1}-\mathbf{U}_{i,j}^{n}}{\Delta t} + \dfrac{(\boldsymbol{\chi}_{i+1/2,j}-\boldsymbol{\chi}_{i-1/2,j})}{\Delta x}+ \dfrac{(\boldsymbol{\chi}_{i,j+1/2}-\boldsymbol{\chi}_{i,j-1/2})}{\Delta y}=\dfrac{1}{2}\left( {\nu_{i,j}^{n}}\mathbf{S}(\mathbf{U}_{i,j}^{n})
    +{\nu_{i,j}^{n+1}}\mathbf{S}(\mathbf{U}_{i,j}^{n+1})\right).
\end{aligned}
\end{equation}

\subsection{Alternative asymptotic diffusion schemes}
In the Maxwellian reconstruction (\ref{reconstruction_maxwellian}), the conservative variables are used to compute the slopes. This particular choice results in a limit scheme (\ref{flux_diffusion}) for to the diffusion equation (\ref{diffusion}) where the energy flux $\boldsymbol\nabla_\mf{x} \frac{q^2}{\rho}$ is expanded in $2\frac{q}{\rho}\boldsymbol\nabla_\mf{x} q-\frac{q^2}{\rho^2}\boldsymbol\nabla_\mf{x} \rho$. Consequently, this procedure leads to a non-natural scheme in which the interface values of the macroscopic quantities ($\rho_{i+1/2,j}^n$ and $q_{i+1/2,j}^n$) occur. Other choices that lead to correct diffusion schemes are possible.

For example, let $\mf{Z}=\begin{pmatrix}\rho & T \end{pmatrix}^T$ be the vector of non-conservative variables. In that scenario, the Maxwellian slopes take the same form as in (\ref{pentes_maxwellienne}) (with finite differences on $\mf{Z}$ in them) and the matrix becomes:
\begin{equation}
    K(\mf{Z})=M[\mf{Z}]\begin{pmatrix}
    \dfrac{1}{\rho} & -\dfrac{3}{2T} \\
    0 & \dfrac{1}{2T^2}
    \end{pmatrix}.
\end{equation}
Using those variables, the limit diffusion fluxes are
\begin{subequations}
    \begin{align}
    \begin{split}
    \Phi^\rho_{i+1/2,j}&\underset{\epsilon \to 0}{\longrightarrow}\dfrac{-1}{\sigma_{i+1/2,j}} \left[T_{i+1/2,j}^n\dfrac{\rho_{i+1,j}^n-\rho_{i,j}^n}{\Delta x}+\rho_{i+1/2,j}^n\dfrac{T_{i+1,j}^n-T_{i,j}^n}{\Delta x}\right],
    \end{split} \\
    \begin{split}
    \Phi^q_{i+1/2,j}&\underset{\epsilon \to 0}{\longrightarrow}\dfrac{-5}{2\sigma_{i+1/2,j}} \left( (T_{i+1/2,j}^n)^2\dfrac{\rho_{i+1,j}^n-\rho_{i,j}^n}{\Delta x} + 2\rho_{i+1/2,j}^n T_{i+1/2,j}^n \dfrac{T_{i+1,j}^n-T_{i,j}^n}{\Delta x}\right),
    \end{split}
    \end{align}
\end{subequations}
which are consistent approximations of the fluxes of the diffusion system (\ref{diffusion}) written in terms of the non-conservative variables:
\begin{equation}
    \begin{cases}
    \partial_t \rho &= \boldsymbol\nabla_{\mf{x}} \left(\dfrac{1}{\sigma}\left( T\boldsymbol\nabla_{\mf{x}}\rho+\rho\boldsymbol\nabla_{\mf{x}}T\right)\right)\\
    \partial_t \dfrac{3}{2}\rho T &= \boldsymbol\nabla_{\mf{x}} \left(\dfrac{5}{2\sigma}\left( T^2\boldsymbol\nabla_{\mf{x}} \rho +2\rho T\boldsymbol\nabla_{\mf{x}} T\right)\right)
    \end{cases}.
\end{equation}
In the two previous schemes, interface values of the macroscopic quantities appear. The natural finite volume scheme involves a direct discretization of the energy flux $\boldsymbol\nabla_\mf{x}\frac{q^2}{\rho}$ without gradient expansion. To obtain this scheme, the set of variable $\mf{Y}=\begin{pmatrix}q & y\end{pmatrix}$, where $y=\frac{q^2}{\rho}$, should be chosen. In that case, the matrix linked to the derivative of the Maxwellian is
\begin{equation}
    K(\mf{Y})=M[\mf{Y}]\begin{pmatrix}
                  \dfrac{7}{2q} & \dfrac{-3}{2y} \\
                  \dfrac{-5}{2y} & \dfrac{3q}{2y^2}
    \end{pmatrix},
\end{equation}
and the limit of the associated macroscopic fluxes are
\begin{subequations}
    \begin{align}
    \begin{split}
    \Phi^\rho_{i+1/2,j}&\underset{\epsilon \to 0}{\longrightarrow}\dfrac{-2}{\sigma_{i+1/2,j}} \dfrac{q_{i+1,j}^n-q_{i,j}^n}{\Delta x},
    \end{split} \\
    \begin{split}
    \Phi^q_{i+1/2,j}&\underset{\epsilon \to 0}{\longrightarrow}\dfrac{5}{3}\dfrac{-2}{3\sigma_{i+1/2,j}} \dfrac{y_{i+1,j}^n-y_{i,j}^n}{\Delta x},
    \end{split}
    \end{align}
\end{subequations}
which is the desired scheme.

\section{An extension to unstructured meshes \label{sec:uns}}
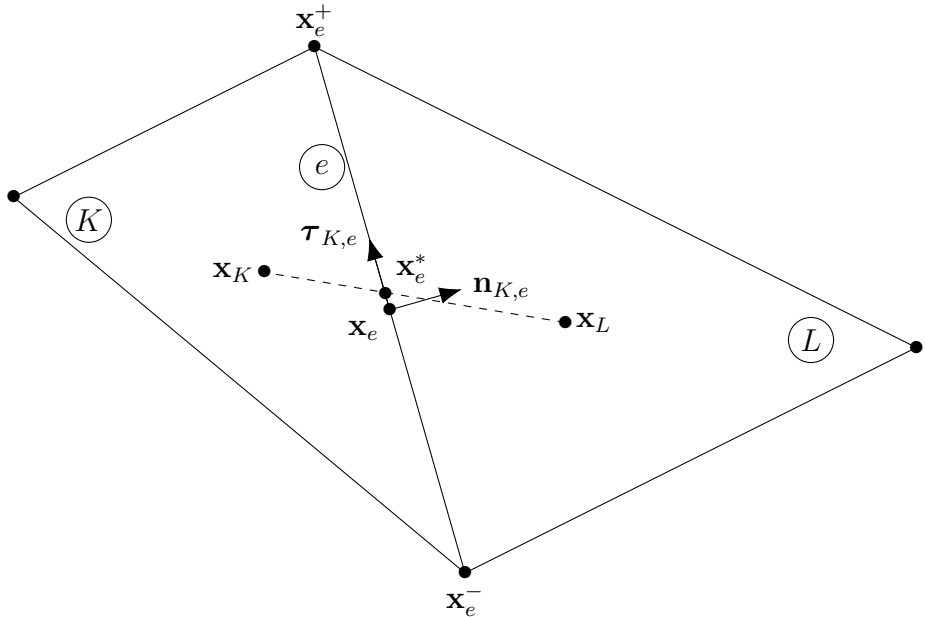
\begin{figure}[htbp!]
    \centering
    \begin{tikzpicture}
        \coordinate (A) at (-6,2);
        \coordinate (B) at (6,0);
        \coordinate (S) at (0,-3);
        \coordinate (N) at (-2,4);
        
        \coordinate (AN) at ($ (A)!.5!(N) $);
        \coordinate (BN) at ($ (B)!.5!(N) $);
        \coordinate (XW) at ($ (S)!.6666!(AN) $);
        \coordinate (XE) at ($ (S)!.6666!(BN) $);
        
        \coordinate (XC) at ($ (S)!.5!(N) $);
        \coordinate (XCP) at ($ (S)!.53!(N) $);
        \coordinate (XCC) at ($ (S)!.53!(N) $);
        \coordinate (XI) at ($ (XE)!.6!(XW) $);
        \coordinate (XEE) at ($ (XW)!2.!(XI) $);
        
        \coordinate (XP) at ($ (S)!(XE)!(N) $);
        \coordinate (XPP) at ($ (S)!(XEE)!(N) $);
        
        \draw (A) node[anchor=east] {} node {$\bullet$};
        \draw (B) node[anchor=west] {} node {$\bullet$};
        \draw (S) node[below] {$\textbf{x}_e^-$} node {$\bullet$};
        \draw (N) node[above] {$\textbf{x}_e^+$} node {$\bullet$};
        \draw (XW) node[anchor=east] {$\textbf{x}_K$} node {$\bullet$};
        \draw (XE) node[anchor=west] {$\textbf{x}_L$} node {$\bullet$};
        \draw (XC) node[below left]{$\textbf{x}_e$} node{$\bullet$};
        \draw (XCP) node[above right]{$\textbf{x}_e^*$} node{$\bullet$};
        
        \draw (S) -- (A) -- (N);
        \draw[dashed] (XE) -- (XW);
        \draw (S) -- (B) -- (N);
        \draw (S) -- (N) node[near end, left]{$$};
        
        \draw[-{Latex[length=3mm]}] (XC) -- ($ (XC)!1cm!270:(N) $) node[anchor=west]{$\mathbf{n}_{K,e}$};
        \draw[-{Latex[length=3mm]}] (XC) -- ($ (XC)!1cm!0:(N) $) node[anchor=east]{$\boldsymbol\tau_{K,e}$};
        
        \draw ($ (A)!0.2!(XC) $) circle [radius=0.3] node {$K$};
        \draw ($ (B)!0.2!(XC) $) circle [radius=0.3] node {$L$};
        \draw ($ (-2.3,4)!0.2!(-0.3,-4) $) circle [radius=0.3] node{$e$};
    \end{tikzpicture}
    \caption{Schematic view of two triangular elements of an unstructured mesh}
    \label{fig:maille_non_struct2}
\end{figure}

In this section an extension to unstructured meshes is proposed for UGKS and therefore UGKS-M1. In the following, we consider a conform (in the sens of finite volume) mesh of the domain $\mathscr{D}$. Moreover, the elements of this mesh can be any polygon but are assumed to be triangles. In our numerical tests, the meshes are generated using GMSH \cite{gmsh}. Notations are given in figure \ref{fig:maille_non_struct2}.

For any non-boundary face $e$, this edge separates two cells, $K$ and $L$. The center of mass of the triangles are denoted by $\mf{x}_K$ and $\mf{x}_L$ and the center of the edge by $\mf{x}_e$. The outer normal of the edge $e$ from triangle $K$ is $\mf{n}_{K,e}$ and the tangential vector is $\boldsymbol\tau_{K,e}$ (the basis $(\mf{n}_{K,e},\boldsymbol\tau_{K,e})$ is direct) and the edge vertices $\mf{x}_e^\pm$ are such that $(\mf{x}_e^+-\mf{x}_e^-)\cdot \boldsymbol\tau_{K,e}>0$. The intersection point between the face and the line formed by the centers of the two cells is denoted by $\mf{x}_e^*$.

\subsection{UGKS}
\subsubsection{Construction of the scheme}
First, the finite volume formulation of the electron transport equation on this mesh is (for any triangle $K$)
\begin{equation}
    \dfrac{f_K^{n+1}-f_K^n}{\Delta t} + \sum_{e\in\mathscr{F}_K} \dfrac{|e|}{|K|} \phi_{K,e} = \nu(\mf{W}_K^{n+1})\left(M[\mf{W}_K^{n+1}]-f_K^{n+1}\right),
\end{equation}
where $f_K^n$ is the mean value of the distribution function in cell $K$ at time $t_n$ and $\phi_{K,e}$ is the mean outward flux from cell $K$ to cell $L$ across the face $e$:
\begin{equation} \label{flux_micro_transv}
    \phi_{K,e} = \dfrac{1}{\eta\Delta t}\int_{t_n}^{t_{n+1}} v \boldsymbol\Omega \cdot \mf{n}_{K,e} f(t,\mf{x}_e,v,\boldsymbol\Omega) \mathrm{d}t. 
\end{equation}
As in the structured case, the macroscopic variables are computed first (in order to evaluate the implicit collision term) using the corresponding macroscopic finite volume formulation:
\begin{equation} \label{vf_macro_ugks_uns}
    \dfrac{\mf{W}_K^{n+1}-\mf{W}_K^n}{\Delta t} + \sum_{e\in\mathscr{F}_K} \dfrac{|e|}{|K|} \boldsymbol\Phi_{K,e} =0,
\end{equation}
where $\boldsymbol\Phi_{K,e}=\ccrochet{\boldsymbol\Psi(v)\phi_{K,e}(v,\boldsymbol\Omega)}$ is the vector of macroscopic fluxes. Contrary to before, the UGKS numerical flux needs to be defined in a arbitrary direction $\mf{n}_{K,e}$. In order to achieve this, the Duhamel formula (\ref{representation_interface}) is first rewritten at the center of the face $\mf{x}_e$. The numerical flux is then obtained by using appropriate reconstructions of the distribution function and of the Maxwellian before time integration. For the sake of simplicity, only the first order version of the scheme is presented here. As a consequence, a constant reconstruction of the distribution function is used:
\begin{equation} \label{reconstruction_f_non_struct}
   \tilde{f}(t_n,\mf{x},\cdot,\cdot)=\begin{cases}
    f_K^n  &\text{ if } \left(\mf{x}-\mf{x}_e\right) \cdot \mf{n}_{K,e} \leq 0\\
    f_L^n  &\text{ if } \left(\mf{x}-\mf{x}_e\right) \cdot \mf{n}_{K,e} \geq 0\\
    \end{cases}.
\end{equation}
The Maxwellian reconstruction is still chosen to be continuous across the face:
 \begin{equation} \label{reconstruction_maxwellian_non_struct}
    \tilde{M}(t,\mf{x},\cdot)=M[\mf{W}_{e}^n]+
    \begin{cases}
    \boldsymbol{\nabla}_\mf{x}^- M_e^n \cdot (\mf{x}-\mf{x}_e) \text{ if } (\mf{x}-\mf{x}_e) \cdot \mf{n}_{K,e} \leq 0 \\
    \boldsymbol{\nabla}_\mf{x}^+ M_e^n \cdot (\mf{x}-\mf{x}_e) \text{ if } (\mf{x}-\mf{x}_e) \cdot \mf{n}_{K,e} \geq 0
    \end{cases}.
\end{equation}
As in the Cartesian case, the gradient $\boldsymbol\nabla_\mf{x}^\pm M_e^n$ can be linked to the gradient of the set of variables chosen to reconstruct the maxwellian. This choice determines the form of the diffusion system that is discretized in the limit. For example, if the conservative variables are used then the gradients are
\begin{equation}\label{slopes}
    \boldsymbol{\nabla}_{\mf{x}}^\pm M_e^n = \begin{pmatrix}\partial_n^\pm M_e^n \\ \partial_\tau^\pm M_e^n \end{pmatrix} = \left(\mf{K}(\mf{W}_e^n) \boldsymbol{\nabla}_\mathbf{x}^\pm \mf{W}_e^n\right)^T \boldsymbol\Psi,
\end{equation}
where the matrix $\mf{K}$ is defined in (\ref{mat_k}). The approximations of the gradient of the conservative variables on both side of the face $\boldsymbol\nabla_\mf{x}^\pm \mf{W}_e^n = \begin{pmatrix}\delta_n^\pm \mf{W}_e^n & \delta_\tau^\pm \mf{W}_e^n \end{pmatrix}$ is discussed in the following sections. Finally, using these two reconstructions, the UGKS microscopic flux is
\begin{equation} \label{ugks_uns}
    \begin{aligned}
        \phi_{K,e}(v,\boldsymbol\Omega) =& A_e^n v \Omega_n \left(f_K^n \mathbb{1}_{\Omega_n > 0}+f_L^n \mathbb{1}_{\Omega_n< 0}\right) \\
        &+ C_e^n v \Omega_n  M[\mf{W}_e^n] \\
        &+ D_e^n v^2 \Omega_n^2 \left(\delta_n^- M_e^n \mathbb{1}_{\Omega_n>0}+\delta_n^+ M_e^n \mathbb{1}_{\Omega_n<0}\right) \\
        &+ D_e^n v^2 \Omega_n\Omega_\tau \left(\delta_\tau^- M_e^n \mathbb{1}_{\Omega_n>0}+\delta_\tau^+ M_e^n \mathbb{1}_{\Omega_n<0}\right),
    \end{aligned}
    \end{equation}
where $\Omega_n=\boldsymbol\Omega \cdot \mf{n}_{K,e}$ and $\Omega_\tau=\boldsymbol\Omega \cdot \boldsymbol\tau_{K,e}$ are the projected directions and where the integration coefficients are interface value of functions (\ref{coeff}) at $\nu_e^n=\frac{1}{2}(\nu_K^n+\nu_L^n)$. Moreover, the two components of the macroscopic flux are
\begin{subequations} \label{flux_macro_non_struct}
    \begin{align}
    \begin{split}
    \Phi^\rho_{K,e}=&A_e^n \ccrochet{v\Omega_n f_k^n \mathbb{1}_{\Omega_n>0} + v\Omega_n f_L^n \mathbb{1}_{\Omega_n<0}} \\
    &+\dfrac{2D_e^n}{3} \dfrac{\delta_n^-q_e^n+\delta_n^+ q_e^n}{2},
    \end{split} \\
    \begin{split}
    \Phi^q_{K,e}&=\dfrac{A_e^n}{2} \ccrochet{v^3\Omega_n f_k^n \mathbb{1}_{\Omega_n>0} + v^3\Omega_n f_L^n \mathbb{1}_{\Omega_n<0}} \\
    &+\dfrac{2D_e^n}{3} \dfrac{5(q_e^n)^2}{3\rho_e^n}\left(
    -\dfrac{\delta_n^-\rho_e^n+\delta_n^+ \rho_e^n}{2\rho_e^n}
    +2\dfrac{\delta_n^-q_e^n+\delta_n^+ q_e^n}{2 q_e^n}
    \right).
    \end{split}
    \end{align}
\end{subequations}

\subsubsection{The general asymptotic limit}
As in the Cartesian case, only the normal part of the Maxwellian slope contributes to the diffusion flux as the tangential part does not appear in the macroscopic fluxes. Without specifying the normal slopes $\partial_n^\pm\mf{W}_e^n$, the asymptotic limit of the macroscopic fluxes are generically
\begin{subequations}
    \begin{align}
    \Phi_{K,e}^\rho &\underset{\epsilon \to 0}{\longrightarrow} -\dfrac{2}{3\sigma_e^n} \dfrac{\delta_n^-q_e+\delta_n^+ q_e}{2}, \\
    \Phi_{K,e}^q &\underset{\epsilon \to 0}{\longrightarrow} -\dfrac{2}{3\sigma_e^n} \dfrac{5(q_e^n)^2}{3\rho_e^n}\left(-\dfrac{\delta_n^-\rho_e+\delta_n^+ \rho_e}{2\rho_e^n}+2\dfrac{\delta_n^-q_e+\delta_n^+ q_e}{2q_e^n}\right).
    \end{align}
\end{subequations}
These expressions indicates that in order to obtain a good diffusion scheme, the approximation of the half slopes must verify a consistency property with the total slope across the interface. More precisely, the mean of the two half slopes should be a consistent approximation of the normal slope between cells K and L:
\begin{equation} \label{prop_consistance}
    \dfrac{1}{2}\left(\delta_n^- \mf{W}_e^n+\delta_n^+ \mf{W}_e^n\right)=\delta_n \mf{W}_e^n,
\end{equation}
where $\delta_n \mf{W}_e^n$ is a given approximation of the normal gradient of the vector of macroscopic variables. The main objective is to correctly define the half slopes to obtain a specific approximation of the normal gradient while verifying the consistency property (\ref{prop_consistance}). In the following sections, several approaches are proposed to do so.

It is worth noting that if a different set of variables is used to reconstruct the Maxwellian (to discretize an alternative form of the diffusion system in the limit), then this consistency property still holds but applies to the newly chosen variables.

\subsubsection{Naive scheme}
To begin with, the first goal is to define the normal half slopes in order to obtain the two point scheme in the diffusion limit. This scheme consists in approximating the normal gradient using a finite difference formula of the gradient between the center of the two cells:
\begin{equation} \label{approx_naive}
    \delta_n\mf{W}_e^n = \dfrac{\mf{W}_L^n-\mf{W}_K^n}{l},
\end{equation}
where $l=||\mf{x}_L-\mf{x}_K||$. This approximation cannot be used in practice as it assumes that the direction between the center of the two cells is collinear with the normal $\mf{n}_{K,e}$, which is not true in general, especially for deformed meshes. Nevertheless, natural definitions of the half slopes are
\begin{subequations} \label{demi_pentes_naives}
    \begin{align}
        \delta_n^-\mf{W}_e^n &= \dfrac{\mf{W}_e^n-\mf{W}_K^n}{l^-}, \\
        \delta_n^+\mf{W}_e^n &= \dfrac{\mf{W}_{L}^n-\mf{W}_e^n}{l^+},
    \end{align}
\end{subequations}
where $l^-=||\mf{x}_e^*-\mf{x}_k||$ and $l^+=||\mf{x}_L-\mf{x}_e^*||$. Thus, the normal slopes of the Maxwellian are
\begin{subequations} \label{demi_pentes_naives_max}
    \begin{align}
        \delta_n^-M_e^n &= \boldsymbol\Psi^T \mf{K}(\mf{W}_e^n)\dfrac{\mf{W}_e^n-\mf{W}_K^n}{l^-}, \\
        \delta_n^+M_e^n &= \boldsymbol\Psi^T \mf{K}(\mf{W}_e^n)\dfrac{\mf{W}_{L}^n-\mf{W}_e^n}{l^+}.
    \end{align}
\end{subequations}
The tangential slopes of the Maxwellian are neglected. By assuming, (\ref{approx_naive})-(\ref{demi_pentes_naives}), it is easy to find that (\ref{prop_consistance}) is satisfied if 
\begin{equation} \label{interafece_modifiee}
    \mf{W}_e^n=-\dfrac{l^+-l^-}{2l}(\mf{W}_L^n-\mf{W}_K^n) + \dfrac{1}{2}\left(\mf{W}_K^n+\mf{W}_L^n\right).
\end{equation}
When $l^+=l^-$, the interface value obtained with this formula coincidences with the usual definition of interface value given by the arithmetic mean. Finally, using definitions (\ref{demi_pentes_naives}) of the half slopes and definition (\ref{interafece_modifiee}) of the interface macroscopic vector, the limits of the macroscopic fluxes are
\begin{subequations}
    \begin{align}
    \Phi_{K,e}^\rho &\underset{\epsilon \to 0}{\longrightarrow} -\dfrac{2}{3\sigma_e^n} \dfrac{q_L^n-q_K^n}{l}, \\
    \Phi_{K,e}^q &\underset{\epsilon \to 0}{\longrightarrow} -\dfrac{2}{3\sigma_e^n} \dfrac{5(q_e^n)^2}{3\rho_e^n}\left(-\dfrac{\rho_L^n-\rho_K^n}{l}+2\dfrac{q_L^n-q_K^n}{l}\right),
    \end{align}
\end{subequations}
which is the desired (bad) scheme.

\subsubsection{Diamond scheme}
In practice, meshes are deformed and the naive scheme can not be employed. However, the ideas presented in the previous section can be reused to construct a more sophisticated scheme. To begin with, we propose a procedure to obtain the diamond scheme as the asymptotic limit of UGKS. 
\paragraph{The diamond scheme for the diffusion system}

\

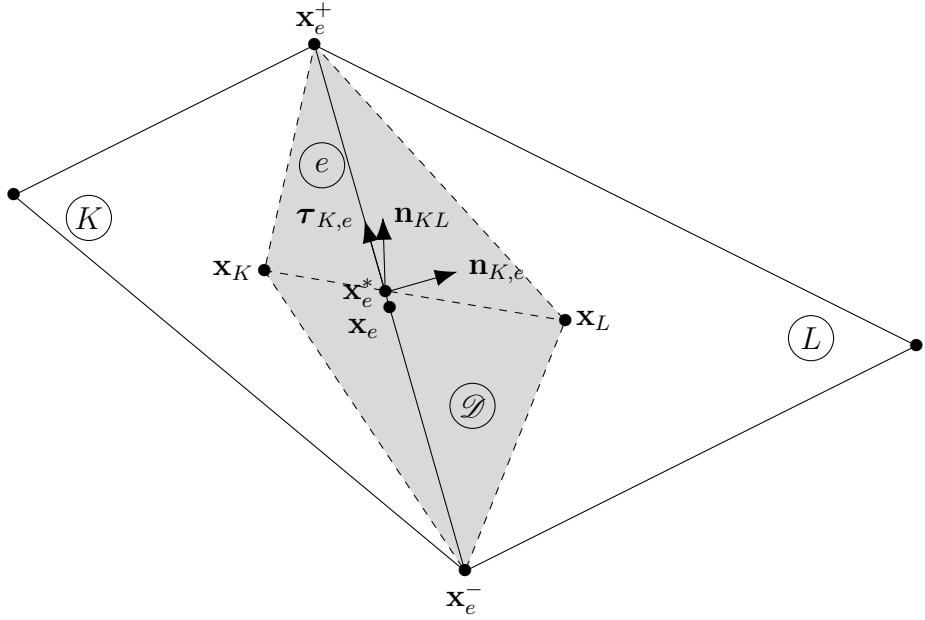
\begin{figure}[ht!]
    \centering
    \begin{tikzpicture}
        \coordinate (A) at (-6,2);
        \coordinate (B) at (6,0);
        \coordinate (S) at (0,-3);
        \coordinate (N) at (-2,4);
        
        \coordinate (AN) at ($ (A)!.5!(N) $);
        \coordinate (BN) at ($ (B)!.5!(N) $);
        \coordinate (XW) at ($ (S)!.6666!(AN) $);
        \coordinate (XE) at ($ (S)!.6666!(BN) $);
        
        \coordinate (XC) at ($ (S)!.5!(N) $);
        \coordinate (XCC) at ($ (S)!.53!(N) $);
        \coordinate (XI) at ($ (XE)!.6!(XW) $);
        \coordinate (XEE) at ($ (XW)!2.!(XI) $);
        
        \coordinate (XP) at ($ (S)!(XE)!(N) $);
        \coordinate (XPP) at ($ (S)!(XEE)!(N) $);
        
        \draw[dashed, fill=gray, fill opacity=0.3] (XW) -- (S) -- (XE) -- (N) -- cycle;
        
        \draw (A) node[anchor=east] {} node {$\bullet$};
        \draw (B) node[anchor=west] {} node {$\bullet$};
        \draw (S) node[below] {$\textbf{x}_e^-$} node {$\bullet$};
        \draw (N) node[above] {$\textbf{x}_e^+$} node {$\bullet$};
        \draw (XW) node[anchor=east] {$\textbf{x}_K$} node {$\bullet$};
        \draw (XE) node[anchor=west] {$\textbf{x}_L$} node {$\bullet$};
        \draw (XC) node[below left]{$\textbf{x}_e$} node{$\bullet$};
        \draw (XCC) node[ left]{$\textbf{x}_e^*$} node{$\bullet$};

        \draw (S) -- (A) -- (N);
        \draw (S) -- (B) -- (N);
        \draw (S) -- (N) node[near end, left]{$$};
        
        \draw[-{Latex[length=3mm]}] (XCC) -- ($ (XCC)!1cm!270:(N) $) node[anchor=west]{$\mathbf{n}_{K,e}$};
        \draw[-{Latex[length=3mm]}] (XCC) -- ($ (XCC)!1cm!0:(N) $) node[anchor=east]{$\boldsymbol\tau_{K,e}$};
         \draw[-{Latex[length=3mm]}] (XCC) -- ($ (XCC)!1cm!282:(XW) $) node[anchor=west]{$\mathbf{n}_{KL}$};
        \draw[dashed] (XW) -- (XE) ;
        
        \draw ($ (A)!0.2!(XC) $) circle [radius=0.3] node {$K$};
        \draw ($ (B)!0.2!(XC) $) circle [radius=0.3] node {$L$};
        \draw ($ (-2.3,4)!0.2!(-0.3,-4) $) circle [radius=0.3] node{$e$};
        \draw (0.1,-0.8) circle [radius=0.3] node{$\mathscr{D}$};
        
    \end{tikzpicture}
    \caption{Schematic view of two triangular elements of an unstructured mesh and the associated diamond $\mathscr{D}=(\mf{x}_K,\mf{x}_e^-,\mf{x}_L,\mf{x}_e^+)$}
    \label{fig:maille_diam}
\end{figure}

This diffusion scheme has first been introduced in \cite{diamant} and consists in approximating the gradient by the mean of well reconstructed gradient in the diamond polyhedron $\mathscr{D}$ (see figure \ref{fig:maille_diam}).
Using linear reconstructions of $\mf{W}$ on the edge of $\mathscr{D}$, the diamond approximation of the normal gradient is defined by:
\begin{equation}
\delta_n\mf{W}_e^n = \dfrac{\mf{W}_L^n-\mf{W}_K^n}{h} - \kappa \dfrac{\mf{W}_e^{n,+}-\mf{W}_e^{n,-}}{e},
\end{equation}
where $h=(\mf{x}_L-\mf{x}_K)\cdot\mf{n}_{K,e}$ is the projected length, $e=||\mf{x}_e^+-\mf{x}_e^-||$ is the edge length, $\kappa=\frac{\mf{x}_L-\mf{x}_K)\cdot\mf{t}_{K,e}}{\mf{x}_L-\mf{x}_K)\cdot\mf{n}_{K,e}}$ is the tangent of an angle and represents the cells deformation and $\mf{W}_e^{n,\pm}$ are the values of the macroscopic variables on the vertices $\mf{x}_e^\pm$. It has been shown in \cite{diamant} that this approximation leads to a second-order accurate scheme. This approximation can then be used in the fluxes of the diffusion system (\ref{diffusion}) to obtain the diamond scheme:
\begin{subequations} \label{limite_diam}
    \begin{align}
    \Phi_{K,e}^\rho &= -\dfrac{2}{3\sigma_e^n} \left(\dfrac{q_L^n-q_k^n}{h} -\kappa \dfrac{q_e^{n,+}-q_e^{n,-}}{|e|}\right), \\
    \Phi_{K,e}^\rho &= -\dfrac{2}{3\sigma_e^n} \dfrac{5(q_e^n)^2}{3\rho_e^n}\left[
    -\dfrac{1}{\rho_e^n}\left(\dfrac{\rho_L^n-\rho_k^n}{h} -\kappa \dfrac{\rho_e^{n,+}-\rho_e^{n,-}}{|e|}\right)
    +\dfrac{2}{q_e^n}\left(\dfrac{q_L^n-q_k^n}{h} -\kappa \dfrac{q_e^{n,+}-q_e^{n,-}}{|e|}\right)
    \right].
    \end{align}
\end{subequations}
It should be noted that if the cells are not deformed, then $\kappa=0$ and the diamond approximation of the normal gradient is the previously introduced two points formula (\ref{approx_naive}). Moreover, as the values of the macroscopic variables at the vertices of the mesh are not primal unknowns of the problem, a reconstruction procedure has to be proposed. Following \cite{diamant} work, a least squared based method is used (see \ref{ann6} for details).

\paragraph{Definition of the UGKS macroscopic states}

\

In order to obtain the diamond scheme as the asymptotic limit of UGKS, the half normal gradients of the macroscopic quantities $\delta_n^\pm\mf{W}_e^n$ in the Maxwellian reconstruction should be defined accordingly. To do this, we propose to follow a similar approach as in the regular diamond scheme construction and to approximate the gradients on both side of the interface by the mean of a well reconstructed gradient on the half diamonds $\mathscr{D}^-=(\mf{x}_K,\mf{x}_e^-,\mf{x}_e^+)$ and $\mathscr{D}^+=(\mf{x}_L,\mf{x}_e^+,\mf{x}_e^-)$. For example, let $\mf{P}^-_e$ be the mean of $\boldsymbol\nabla_\mf{x} \mf{W}$ in $\mathscr{D}^-$ at time $t_n$. Using linear reconstructions of $\mf{W}$ on the edges of $\mathscr{D}^-$, it can be shown that $\mathbf{P}_e^-$ is the unique solution of the system:
\begin{equation*}
\begin{cases}
    \mathbf{P}_e^- (\mathbf{x}_e^*-\mathbf{x}_K) &= \mf{W}_e^n-\mf{W}_K^n \\
    \mathbf{P}_e^- (\mathbf{x}_e^+-\mathbf{x}_e^-) &= \mf{W}_e^{n,+}-\mf{W}_e^{n,-} \\
\end{cases},
\end{equation*}
which is
\begin{equation}
    \mathbf{P}^-_e = \left(\dfrac{\mf{W}_e^n-\mf{W}_K^n}{h^-}-\kappa \dfrac{\mf{W}_e^{n,+}-\mf{W}_e^{n,-}}{|e|}\right) \otimes \mathbf{n}_{K,e} + \dfrac{\mf{W}_e^{n,+}-\mf{W}_e^{n,-}}{|e|} \otimes \mathbf{t}_{K,e},
\end{equation}
where $h^-=(\mf{x}_e^*-\mf{x}_K)\cdot \mf{n}_{K,e}$ is the half projected distance. Thus, in UGKS, the normal slopes of the Maxwellian are
\begin{subequations} \label{demi_pente_diam}
    \begin{align}
        \delta_n^-M_e^n &=\boldsymbol\Psi^T\mf{K}(\mf{W}_e^n)\left(\dfrac{\mf{W}_e^n-\mf{W}_K^n}{h^-}-\kappa \dfrac{\mf{W}_e^{n,+}-\mf{W}_e^{n,-}}{|e|}\right), \\
        \delta_n^+M_e^n &=\boldsymbol\Psi^T\mf{K}(\mf{W}_e^n)\left(\dfrac{\mf{W}_L^n-\mf{W}_e^n}{h^+}-\kappa \dfrac{\mf{W}_e^{n,+}-\mf{W}_e^{n,-}}{|e|}\right),
    \end{align}
\end{subequations}
where $h^+=(\mf{x}_L-\mf{x}_e^*)\cdot \mf{n}_{K,e}$. The tangential part of $\mf{P}_e^-$ can also be employed in the Maxwellian reconstructions but is neglected in practice. Using these definitions of the half slopes, the consistency property (\ref{prop_consistance}) is not ensured in general (as $h^+$ is \textit{a priori} different than $h^-$). As in the previous section, this problem can be addressed by defining the interface macroscopic state $\mf{W}_e^n$ correctly. After some algebra, it can be shown that the following definition:
\begin{equation} \label{interafece_modifiee_bis}
    \mf{W}_e^n=-\dfrac{h^+-h^-}{2h}(\mf{W}_L^n-\mf{W}_K^n) + \dfrac{1}{2}\left(\mf{W}_K^n+\mf{W}_L^n\right),
\end{equation}
allows to ensure the consistency property. Finally, the limits of the UGKS macroscopic fluxes are the excepted diamond scheme (\ref{limite_diam}) for the diffusion system (\ref{diffusion}).

\paragraph{Modified diamond scheme}

\

In the previous section, the elaboration of the limit scheme requires to define the interface value correctly in order to ensure the consistency property (\ref{prop_consistance}). If other definitions are adopted (such as the kinetic one: $\mf{W}_e^n=\ccrochet{\boldsymbol\Psi \left(f_K^n \mathbb{1}_{\Omega_n>0}+f_L^n \mathbb{1}_{\Omega_n<0}\right)}$), an alternative approach needs to be proposed. As in the Cartesian case, if the mesh has the right geometric property then the consistency property is ensured for any definition of $\mf{W}_e^n$. This occurs when the two cell centers are equidistant to the face ($h^-=h^+$). Thus, a natural idea is to introduce a new point $\mf{x}_{K}^*$ or $\mf{x}_L^*$ which is the symmetric of the nearest center to the face with respect to $\mf{x}_e^*$ (see Figure \ref{fig:diam_modif}). As the distance from this point to $\mf{x}_e^*$ is $l_0=\min(l^-,l^+)$, it allows to define a modified diamond $\mathscr{D}^*$ with enough regularity to ensure the consistency of the slopes.
\begin{figure}[ht!]
    \centering
    \begin{tikzpicture}
        \coordinate (A) at (-6,0);
        \coordinate (B) at (8,-6);
        \coordinate (S) at (0,-4);
        \coordinate (N) at (-2,4);
        
        \coordinate (AN) at ($ (A)!.5!(N) $);
        \coordinate (BN) at ($ (B)!.5!(N) $);
        \coordinate (XW) at ($ (S)!.6666!(AN) $);
        \coordinate (XE) at ($ (S)!.6666!(BN) $);
        
        \coordinate (XC) at ($ (S)!.5!(N) $);
        \coordinate (XI) at ($ (XE)!.6!(XW) $);
        \coordinate (XEE) at ($ (XW)!2.!(XI) $);
        
        \coordinate (XP) at ($ (S)!(XE)!(N) $);
        \coordinate (XPP) at ($ (S)!(XEE)!(N) $);
        
        \draw[dashed, fill=gray, fill opacity=0.3] (XW) -- (S) -- (XE) -- (N) -- cycle;
        \draw[dashed, fill=gray, fill opacity=0.5] (XW) -- (S) -- (XEE) -- (N) -- cycle;
        
        \draw (A) node[anchor=east] {$\mathbf{x}_K^-$} node {$\bullet$};
        \draw (B) node[anchor=west] {$\mathbf{x}_L^+$} node {$\bullet$};
        \draw (S) node[below] {$\textbf{x}_e^-$} node {$\bullet$};
        \draw (N) node[above] {$\textbf{x}_e^+$} node {$\bullet$};
        \draw (XW) node[anchor=east] {$\textbf{x}_K=\mf{x}_K^*$} node {$\bullet$};
        \draw (XE) node[anchor=west] {$\textbf{x}_L$} node {$\bullet$};
        \draw (XC) node[anchor=east] {$\textbf{x}_e$} node {$\bullet$};
        \draw (XI) node[below left] {$\textbf{x}_e^*$} node {$\bullet$};
        \draw (XEE) node[anchor=south] {$\textbf{x}_L^*$} node {$\bullet$};
        
        \draw (S) -- (A) -- (N);
        \draw (S) -- (B) -- (N);
        \draw (XE) -- (XW);
        \draw (S) -- (N) node[near end, left]{$e$};
        
        \draw[-{Latex[length=3mm]}] (XC) -- ($ (XC)!1cm!270:(N) $) node[anchor=west]{$\mathbf{n}_{K,e}$};
        \draw[-{Latex[length=3mm]}] (XI) -- ($ (XI)!1cm!90:(XE) $) node[below right, near end]{$\mathbf{n}_{K,WE}$};
        \draw[-{Latex[length=3mm]}] (XC) -- ($ (XC)!1cm!0:(N) $) node[anchor=west]{$\boldsymbol\tau_{K,e}$};
        

        
    
        \draw ($ (A)!0.2!(XC) $) circle [radius=0.3] node {$K$};
        \draw ($ (B)!0.2!(XC) $) circle [radius=0.3] node {$L$};
        \draw (1.3,-2.1) circle [radius=0.3] node {$\mathscr{D}$};
        \draw (0.2,-2) circle [radius=0.3] node {$\mathscr{D}^*$};
        
    \end{tikzpicture}
    \caption{Schematic view of the diamond and of the modified diamond in a unstructured mesh}
    \label{fig:diam_modif}
\end{figure}
Let $\mf{x}_{K}^*$ and $\mf{x}_L^*$ be the new vertices of the diamond after this modification (in figure \ref{fig:diam_modif}, $\mf{x}_L^*$ is such that $||\mf{x}_L^*-\mf{x}_e^*||=||\mf{x}_e^*-\mf{x}_K||$ and $\mf{x}_K^*=\mf{x}_K$) and let $\mf{W}_{K}^{n,*}$ and $\mf{W}_{L}^{n,*}$ be their corresponding value. The reconstructed gradients on the half diamonds are
\begin{subequations} \label{gradient_demi_diamant_modif}
\begin{align}
        \mathbf{P}^-_{e} &= \left(\dfrac{\mf{W}_e^n-\mf{W}_K^{n,*}}{h_0}-\kappa \dfrac{\mf{W}_e^{n,+}-\mf{W}_e^{n,-}}{|e|}\right) \otimes \mathbf{n}_{K,e} + \dfrac{\mf{W}_e^{n,+}-\mf{W}_e^{n,-}}{|e|} \otimes \mathbf{t}_{K,e}, \\
        \mathbf{P}^+_{e} &= \left(\dfrac{\mf{W}_L^{n,*}-\mf{W}_e^n}{h_0}-\kappa \dfrac{\mf{W}_e^{n,+}-\mf{W}_e^{n,-}}{|e|}\right) \otimes \mathbf{n}_{K,e} + \dfrac{\mf{W}_e^{n,+}-\mf{W}_e^{n,-}}{|e|} \otimes \mathbf{t}_{K,e},
\end{align}
\end{subequations}
where $h_0=\min(h^+,h^-)$ is the minimum of the projected lengths. The half slopes of the Maxwellian are still defined according to (\ref{slopes}). In that case, the consistency property is obtained for any value of the interface value. Finally, the new asymptotic limit of the macroscopic fluxes are
\begin{subequations} \label{limite_diam_modif}
    \begin{align}
    \Phi_{K,e}^\rho &\underset{\epsilon \to 0}{\longrightarrow} -\dfrac{2}{3\sigma_e^n} \left(\dfrac{q_L^{n,*}-q_K^{n,*}}{2h_0} -\kappa \dfrac{q_e^{n,+}-q_e^{n,-}}{|e|}\right), \\
    \Phi_{K,e}^\rho &\underset{\epsilon \to 0}{\longrightarrow} -\dfrac{2}{3\sigma_e^n} \dfrac{5(q_e^n)^2}{3\rho_e^n}\left[
    -\dfrac{1}{\rho_e^n}\left(\dfrac{\rho_L^{n,*}-\rho_K^{n,*}}{2h_0} -\kappa \dfrac{\rho_e^{n,+}-\rho_e^{n,-}}{|e|}\right)
    +\dfrac{2}{q_e^n}\left(\dfrac{q_L^{n,*}-q_K^{n,*}}{2h_0} -\kappa \dfrac{q_e^{n,+}-q_e^{n,-}}{|e|}\right)
    \right].
    \end{align}
\end{subequations}
It should be noted that the limit scheme is no longer the diamond scheme due to the new point $\mf{W}_{KL}^*$. Moreover, in order to use this scheme, it is necessary to compute the value of the macroscopic variables at this point, as a function of the problem unknowns. Each component $u_K^*$ of $\mf{W}_K^*$ is defined by $u^*_K=\tilde{u}(\mf{x}_K^*)$ with $\tilde{u}(\mf{x})=u_K + \mathbf{a}\cdot (\mathbf{x}-\mathbf{x}_K)$, where the vector $\mf{a}$ is defined so as to minimise the quadratic error $E=\frac{1}{2}\sum_{i\in\mathscr{N}_K} |\tilde{u}(\mf{x}_{K})-u_i|^2$, and where the nodal values $u_i$ are themselves defined by the least square method as in the previous diamond scheme.

In our tests, no significant differences could be highlighted between the modified diamond scheme and the normal one. Moreover, the influence of the definition of the macroscopic vector at the interface also seems to be negligible. As a consequence, and for practical reasons, this approach will not be further studied.

\subsection{UGKS-M1}
Now that UGKS has been established on unstructured meshes, the implementation of UGKS-M1 is straightforward. First, the new microscopic finite volume formulation is
\begin{equation}
    \dfrac{\mf{U}_K^{n+1}-\mf{U}_K^n}{\Delta t} + \sum_{e \in \mathscr{F}_K} \dfrac{|e|}{|K|}\boldsymbol\chi_{K,e}= \nu(\mf{W}_K^{n+1}) \mf{S}(\mf{U}_k^{n+1}),
\end{equation}
where $\mf{U}_K^n$ is the mean of the M1 variables in cell $K$ at time $t_n$, $\boldsymbol\chi_{K,e}=\crochet{\mf{m}\phi_{K,e}}$ is the mesoscopic fluxes vector in direction $\mf{n}_{K,e}$ and $\mf{S}(\mf{U}_k^{n+1})=\begin{pmatrix}M_0[\mf{W}_{i,j}^{n+1}] & 0 \end{pmatrix}^T$ is the source term. The UGKS macroscopic finite volume formulation (\ref{vf_macro_ugks_uns}) is used to compute the macroscopic variables $\mf{W}_{i,j}^{n+1}$ at time $t_{n+1}$. The mesoscopic numerical fluxes are still moments of the UGKS microscopic flux (\ref{ugks_uns}) computed using the closed distribution function $\hat{f}_k^n=\hat{f}(\mf{U}_k^n)$:
\begin{subequations} \label{te_flux_micro1_ugksm1_unstr}
    \begin{align}
    \begin{split}
    \chi^0_{i+1/2,j}(v)=&A_{e}^n v \crochet{\Omega_n \hat{f}_K^{n} \mathbb{1}_{\Omega_n> 0}+\Omega_n\hat{f}_L^{n} \mathbb{1}_{\Omega_n< 0}} \\
    +&\dfrac{D_{e}^n}{6} v^2 \mf{K}_0\left(\mf{W}_{e}^n\right) \left(\delta_n^-\mf{W}_e^n + \delta_n^+\mf{W}_e^n \right) \cdot \boldsymbol\Psi,
    \end{split} \\
    \begin{split}
    \boldsymbol\chi^1_{i+1/2,j}(v)=&A_{e}^n v \crochet{\Omega_n\boldsymbol\Omega \hat{f}_{K}^{n} \mathbb{1}_{\Omega_n> 0}+\Omega_n\boldsymbol\Omega\hat{f}_{L}^{n} \mathbb{1}_{\Omega_n< 0}} \\
    +&\dfrac{C_e^n}{3} v M_0[\mf{W}_{e}^n] \mathbf{n}_{K,e} \\
    +&\dfrac{D_e^n}{8} v^2 \mf{K}_0\left(\mf{W}_{e}^n\right) \left(\delta_n^-\mf{W}_e^n - \delta_n^+\mf{W}_e^n \right) \cdot \boldsymbol\Psi \mf{n}_{K,e}.
    \end{split}
\end{align}
\end{subequations}
Moreover, the macroscopic variables are previously updated with the finite volume formulation (\ref{vf_macro_ugks_uns}), where the fluxes are identical to the UGKS ones (\ref{flux_macro_non_struct}) but with a closed distribution function:
\begin{subequations} \label{flux_macro_non_struct_M1}
    \begin{align}
    \begin{split}
    \Phi^\rho_{K,e}=&A_e^n \ccrochet{v\Omega_n \hat{f}_k^n \mathbb{1}_{\Omega_n>0} + v\Omega_n \hat{f}_L^n \mathbb{1}_{\Omega_n<0}} \\
    &+\dfrac{2D_e^n}{3} \dfrac{\delta_n^-q_e^n+\delta_n^+ q_e^n}{2},
    \end{split} \\
    \begin{split}
    \Phi^q_{K,e}&=\dfrac{A_e^n}{2} \ccrochet{v^3\Omega_n \hat{f}_k^n \mathbb{1}_{\Omega_n>0} + v^3\Omega_n \hat{f}_L^n \mathbb{1}_{\Omega_n<0}} \\
    &+\dfrac{2D_e^n}{3} \dfrac{5(q_e^n)^2}{3\rho_e^n}\left(
    -\dfrac{\delta_n^-\rho_e^n+\delta_n^+ \rho_e^n}{2\rho_e^n}
    +2\dfrac{\delta_n^-q_e^n+\delta_n^+ q_e^n}{2 q_e^n}
    \right).
    \end{split}
    \end{align}
\end{subequations}
The half slopes of the Maxwellian are defined according to the desired limit scheme. Finally, the half moments of the M1 distribution function are computed using the method presented in section \ref{sec:quad}.

\section{Numerical results \label{sec:results}}
\subsection{General framework of the test cases}
In the following section, numerical test cases are presented to evaluate UGKS-M1 on structured and unstructured meshes in comparaison with a standard HLL scheme. Several test cases are selected to validate the properties of the scheme and highlight its capabilities. In most test cases, the initial state is at equilibrium: the distribution functions are Maxwellians at a given macroscopic state $\mathbf{W}_{i,j}^0$ and the corresponding M1 variables are
\begin{equation}
    \mathbf{U}_{i,j}^0=\begin{pmatrix}
    \crochet{M[\mathbf{W}_{i,j}^0]} \\
    \crochet{\boldsymbol\Omega M[\mathbf{W}_{i,j}^0]}
    \end{pmatrix}=\begin{pmatrix}
    4\pi M[\mathbf{W}_{i,j}^0] \\
    0 
    \end{pmatrix}.
\end{equation}
The second variable $\mathbf{f}_1$ can be set to $f_0u  \mathbf{d}$ to simulate a non equilibrium initial situation and to induce particle transport in direction $\mathbf{d} \in \mathbb{R}^3$ with an anisotropic factor $u<1$. The simulations parameters are summarised in table \ref{tab:param} where $H$ is a regularised Heaviside step function defined as
\begin{equation}
    \forall x\in\mathbb{R},~H(x,x_0,T_0,T_1,e)=T_0 + \dfrac{T_1-T_0}{2}\left(\dfrac{2}{\pi}\arctan\left(\dfrac{x-x_0}{e}\right)+1\right).
\end{equation}
The temperature profile used for the non-local test cases is
\begin{equation}
    \forall x\in [0,1],~T^{NL}(x)=H(x,0.5,1,4,0.01) \left( \mathbb{1}_{x<0.5} +\left(\dfrac{0.05}{(x+0.5)^6}+0.95\right) \mathbb{1}_{x>0.5}
    \right).
\end{equation}

\begin{table}[!ht]
    \center
    \begin{tabular}{| c || l | l | l | l |}
     \hline			
        ~ & $\eta$ & $\epsilon$ & \makecell[l]{Initial \\ Condition} & \makecell[l]{Boundary \\ Conditions} \\\hline
       Convergence & $1$ & $1$ & $
       \begin{matrix}
       \rho(\textbf{x})=&1 \\ T(\textbf{x}) =& 1 + 0.25\sin{(2\pi x)} \\ {\mathbf{f}_1}(\mathbf{x},v)=&0.5f_0 \mathbf{e}_x \end{matrix}$ & Periodic\\\hline
       Kinetic  & $1$ & $1$ & $
       \begin{matrix}
       \rho(\textbf{x})=&1 \\ T(\textbf{x}) =&H(x,0.5,1,2,0.001) \\ \mf{f}_1(\mf{x},v)=&0 \end{matrix}$ & Neumann\\\hline
       Intermediate  & $10^{-2}$ & $10^{-2}$ & $
       \begin{matrix}
       \rho(\textbf{x})=&1 \\ T(\textbf{x}) =&H(x,0.5,1,2,0.001) \\ \mf{f}_1(\mf{x},v)=&0 \end{matrix}$ & Neumann\\\hline
       Diffusion  & $10^{-8}$ & $10^{-8}$ & $
       \begin{matrix}
       \rho(\textbf{x})=&1 \\ T(\textbf{x}) =&H(x,0.5,1,2,0.001) \\ \mf{f}_1(\mf{x},v)=&0 \end{matrix}$ & Neumann\\\hline
       \makecell{Energy \\ deposition}  & $1$ & $1$ & $
       \begin{matrix}
       \rho(\textbf{x})=&H(||\mf{x}-\mf{1}||^2,\frac{\sqrt{2}}{2},5,100,0.01) \\ T(\textbf{x}) =&0.5 \\ \mf{f}_1(\mf{x},v)=&0 \end{matrix}$ & \makecell[l]{Neumann: N,E \\ Dirichlet: S,W}\\\hline
        Non-local 1D  & $1$ & $1$ & $
       \begin{matrix}
       \rho(\textbf{x})=&1 \\ T(\textbf{x}) =& T^{NL}(x) \\ \mf{f}_1(\mf{x},v)=&0 \end{matrix}$ & Neumann\\\hline
       Non-local 2D  & $1$ & $1$ & $
       \begin{matrix}
       \rho(\textbf{x})=&1 \\ T(\textbf{x}) =& T^{NL}(x) e^{-y^4} \\ \mf{f}_1(\mf{x},v)=&0 \end{matrix}$ & Neumann\\\hline
     \end{tabular}
     \caption{Simulation parameters\label{tab:param}}
\end{table}
The spatial domain is $\mathscr{D}=[0,1]\times[0,1]$. On structured meshes and for the mono-dimensional test cases, $200$ points in space are used in the horizontal direction and $4$ points in the vertical one. For the bi-dimensional test cases, $50$ points in space are used in both directions. The unstructured meshes are generated with GMSH \cite{gmsh} and the cells are triangles. In order to study the scheme behaviour, two kinds of meshes are considered: non deformed and deformed ones (see figures \ref{fig:mesh_ligne} and \ref{fig:mesh_def} for the meshes used for the mono-dimensional test cases). On non deformed meshes, the cells are close to equilateral triangles. In general, to compare the solutions on structured meshes to those one unstructured meshes, the number of triangles is adjusted to ensure the same cell density.

In every simulation, the half sphere Gauss-Legendre quadrature is performed on $10$ points and the velocity quadrature on $50$ points with a maximum velocity of $12$. The Maxwellian slopes are defined in order to obtain, in the diffusion limit, the diamond scheme for the diffusion system written in developed form. Moreover, we assume that the numerical scheme remain stable under a CFL-like condition combining both a hyperbolic and parabolic condition \cite{mieussens2013}:
\begin{equation}
    \Delta t \leq CFL \left(\dfrac{\min{(\Delta x,\Delta y)}}{v_m} + \dfrac{9}{10}\sigma_0\min{(\Delta x,\Delta y)}^2 \right),
\end{equation}
where $\sigma_0$ is the assumed smallest opacity in the domain at all times. The CFL number is fixed at $0.3$.

\clearpage

\subsection{Test n°1: Relaxation of a sinusoid in a infinite domain}
In this 1D test case, a convergence study is performed on UGKS-M1 on the various discretization steps with a structured mesh in space. A 1D regular initial condition is considered in a bounded domain with periodic boundary conditions. The initial macroscopic state is
a sinusoidal temperature wave and a constant density. The microscopic state is non-Maxwellian to induce particle transport in a preferred direction. 

\begin{figure}
    \centering\input{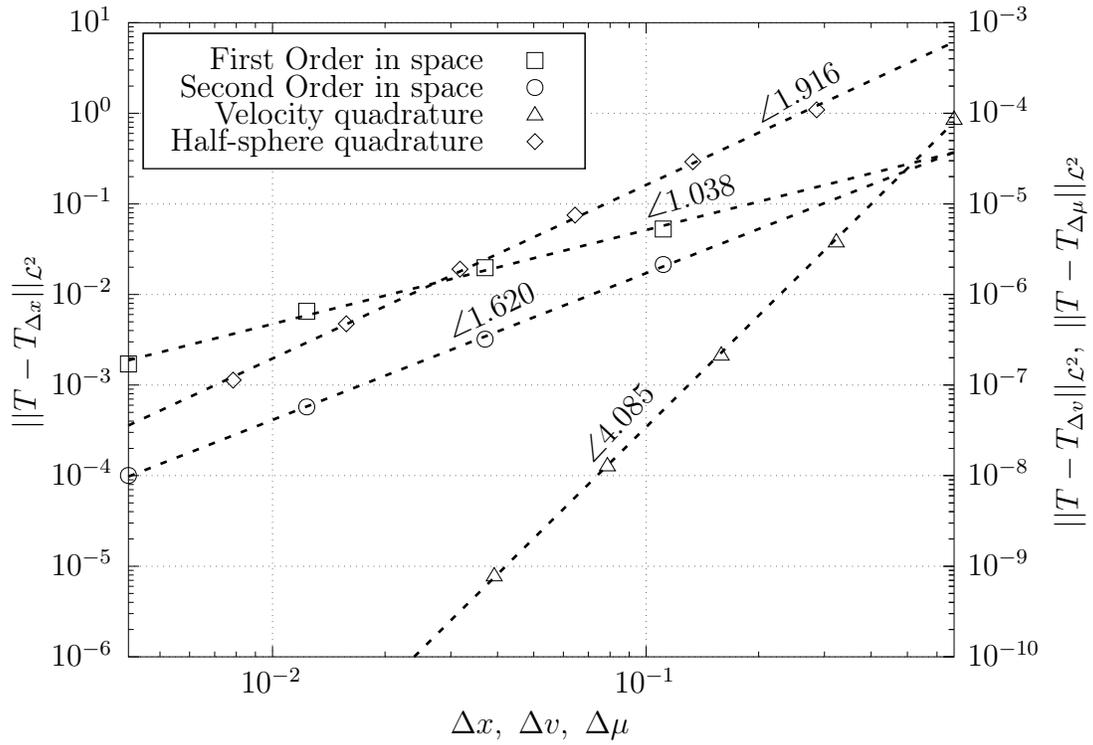}
    \caption{Test n°1: Convergence. UGKS-M1 temperature error as a function of the space step, the velocity step and the half-sphere quadrature step\label{fig:ordre}}
\end{figure}

In figure \ref{fig:ordre}, the $\mathrm{L}^2$ norm of the temperature error at time $t=1$ is represented as a function of the space step, the velocity step and the half-sphere quadrature step. The error is computed based on a reference solution that is sufficiently converged to assume that its difference with respect to the exact solution is negligible. For the second order scheme, the Van Leer limiter is used \cite{vanleer}. Multiple linear regressions are performed to evaluate the convergence orders.

For the space convergence study, we can notice that the fist order UGKS-M1 scheme has the correct convergence rate. With the second order extension, the rate is only $1.6$ despite the regular initial condition. This phenomenon can be explained by the treatment of high velocities that may become numerically stiff due to the lack of particles. Indeed, the Jacobian matrix (\ref{jacob}) in the second order terms is not well defined when $f_0$ is too low, which occurs when the speed $v$ is high. As a consequence, a numerical threshold on $f_0$ ($\epsilon=10^{-8}$) must be introduced to remove this term as it becomes stiff. This procedure degrades the space convergence order.

The velocity convergence order is $4$ and the one on the half-sphere quadrature is $1.9$. The corresponding errors are lower compared to the space one. This demonstrate that much less points are necessary for both the quadrature and the velocity discretization than for the space one to obtain errors of the same order of magnitude. 

\subsection{Test n°2: Transport regime with Neumann boundary conditions}
\begin{figure}[ht!]
    \centering\input{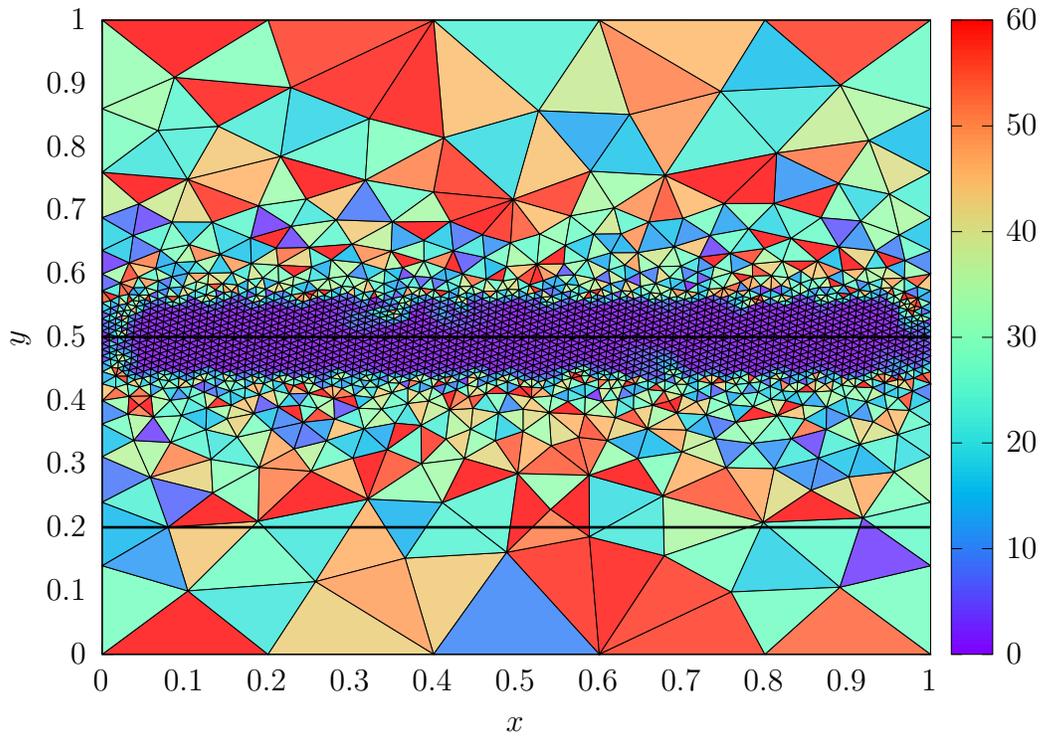}
    \caption{Regular unstructured mesh (used for the 1D test cases) with 4036 triangles and level of angular deformation of triangles in colour scale}
    \label{fig:mesh_ligne}
\end{figure}

\begin{figure}[ht!]
    \centering \input{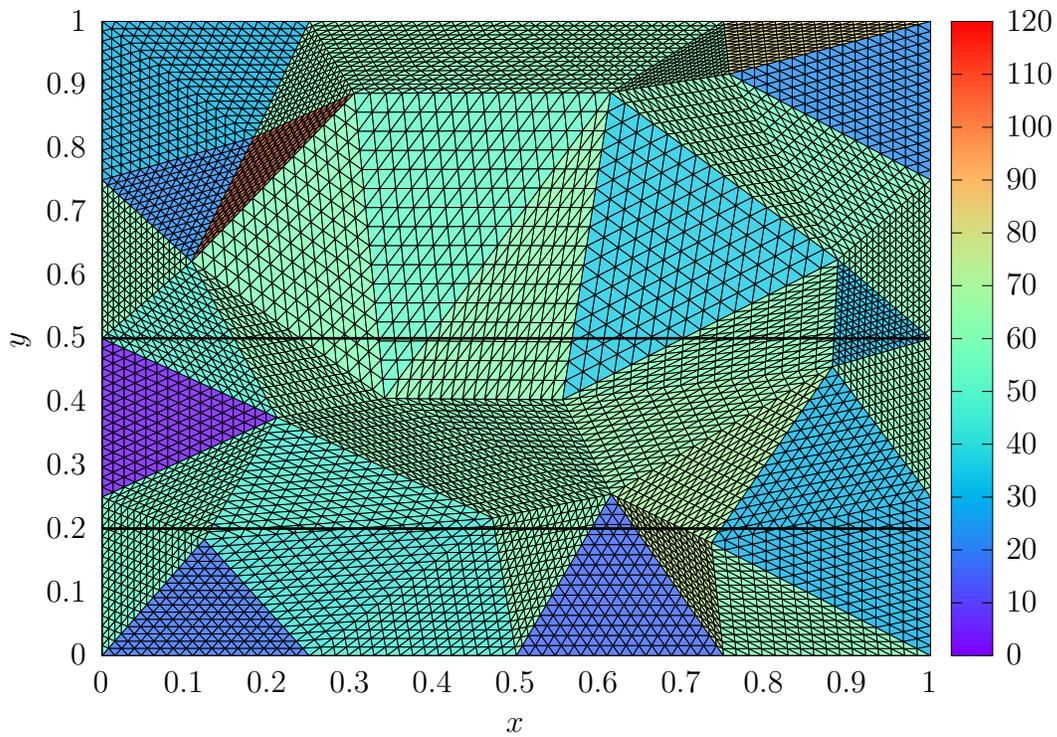}
    \caption{Deformed unstructured mesh and cut lines (used for the 1D test cases) with 10752 triangles and level of angular deformation of triangles in colour scale}
    \label{fig:mesh_def}
\end{figure}

In this 1D test case, we study the relaxation of a temperature step with Neumann boundary conditions in the kinetic regime. In figure \ref{fig:1dkinetic}, the temperature and the density solutions given by UGKS-M1 are represented at different times on a structured mesh and on two unstructured meshes (a regular and a deformed one as shown in figures \ref{fig:mesh_ligne} and \ref{fig:mesh_def}). For the deformed unstructured mesh, the numerical solution is show allong two different cut lines (see figure \ref{fig:mesh_def}), denoted as C1 and C2. First, we can notice that the solution on the structured mesh and on the regular unstructured mesh are almost identical. However, a significant gap can be observed on the deformed mesh. Indeed, the density wave is shifted and the maximum is reduced on both cuts. Additionally, the solution is not constant along vertical lines in the domain despite the symmetry of the problem. In our tests, we were able to show that the error decreases using finer (still deformed) meshes. This shows that this phenomenon is due to numerical diffusion which is induced by the nature of the transport scheme in UGKS-M1 on unstructured mesh. In order to improve these results, we believe that a better distribution function reconstruction should be employed in UGKS to take into account a tangential contribution.

UGKS-M1 has also been compared to a standard HLL scheme and showed similar results with overall less numerical diffusion.

\subsection{Test n°3: Intermediate regime with Neumann boundary conditions}
The same test case is performed in a intermediate regime with a Knudsen number equal to $10^{-2}$. In figure \ref{fig:1dinter2}, the density and temperatures solutions show that the UGKS-M1 on the structured and on the regular unstructured meshes are identical. However, the same phenomenon occurs on the deformed mesh. Even in an intermediate regime, the numerical diffusion of the transport scheme is still clearly visible and affects the quality of the solution.

\subsection{Test n°4: Diffusion regime with Neumann boundary conditions}
Finally, the same test case is performed in the diffusion limit. In figures \ref{fig:diffusion1} and \ref{fig:diffusion2}, two different versions of UGKS-M1 are compared. In the first one, the Maxwellian slopes are defined in order to obtain the naive scheme in the diffusion limit. In the other one, the slopes are set in order to get the diamond scheme. In both case, the solutions on the structured mesh and the regular unstructured one are identical. However, we can clearly notice that the naive scheme is very inaccurate for deformed meshes. Unlike in previous test cases, the observed gap is not due to numerical diffusion but rather results from the inconsistency of the limit scheme when the cells are deformed. The diamond scheme allows to completely correct this problem as the solutions on the deformed meshes are almost identical to the others. Some numerical diffusion is induced but vanishes with mesh refinement.

UGKS-M1 has also been compared with the standard scheme for the diffusion system, and no differences could be highlighted. Moreover, the influence of the choice of variables used to reconstruct the Maxwellian (which determines the form of the discretized system) seems to be completely negligible, at least for these regular initial conditions.

\subsection{Test n°5: Energy deposition in a high density region}
\begin{figure}
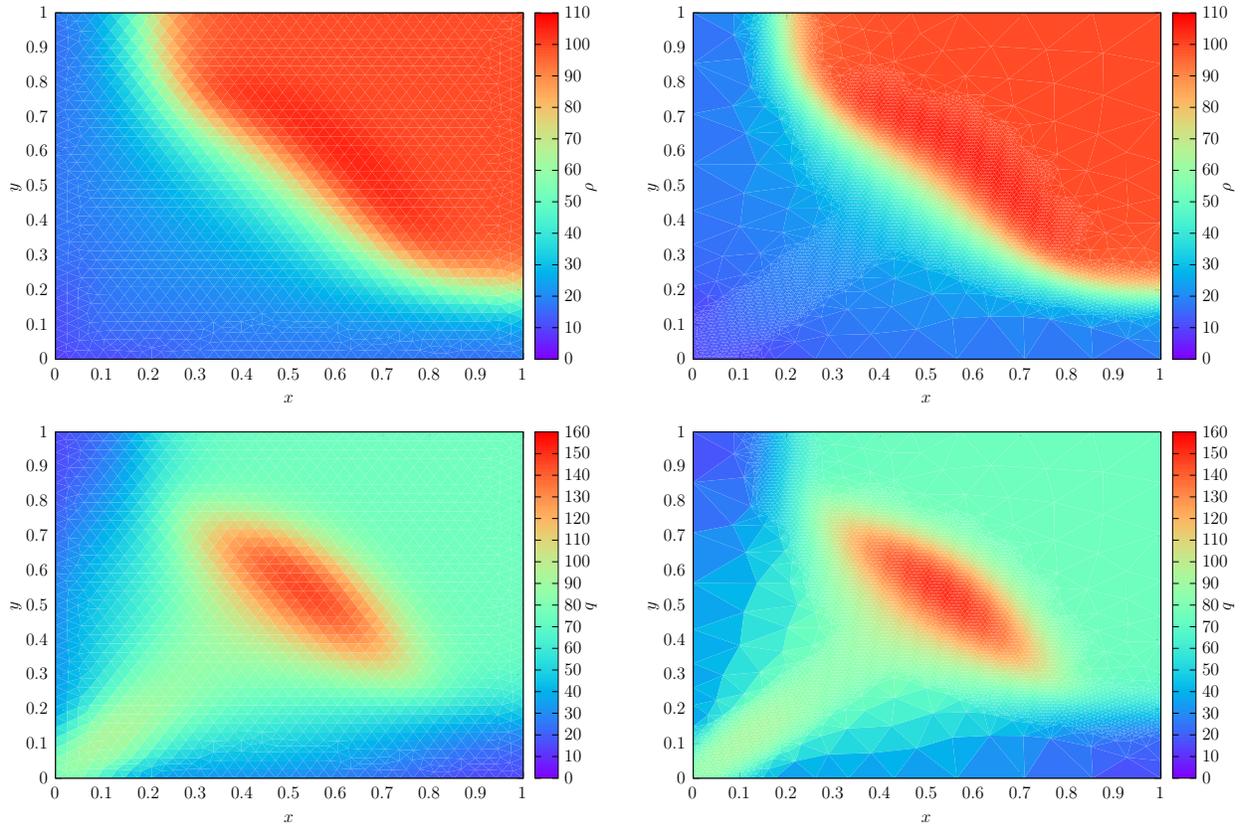

    \centering\scalebox{0.55}{\begin{tabular}{c c}
         \input{density_regular} & \input{density_adapted} \\
         \input{energy_regular} & \input{energy_adapted} \\
    \end{tabular}}
    \caption{Test n°5: Energy deposition in a high density region. Density (above) and energy (below) in the domain for the first order UGKS-M1 on a regular unstructured mesh (left, 4184 triangles) and an adapted one (right, 9726 triangles).\label{fig:uns}}
\end{figure}

In this test case, the scheme capabilities are demonstrated in a 2D situation. At the initial state, the domain is composed of a low and a high density region at a low uniform temperature. On the western and southern boundaries, a M1 distribution function is imposed:
\begin{equation}
    \mf{U}_{S}^{n}(x)=\begin{pmatrix}
        4\pi M[\mf{W}_b(x)] \\
        4\pi M[\mf{W}_b(x)] u \mf{d} 
    \end{pmatrix},\quad
    \mf{U}_{W}^{n}(y)=\begin{pmatrix}
        4\pi M[\mf{W}_b(y)] \\
        4\pi M[\mf{W}_b(y)] u \mf{d} 
    \end{pmatrix},
\end{equation}
where $u=0.95$, $\mf{d}=\dfrac{\sqrt{2}}{{2}}(\mathbf{e}_x+\mathbf{e}_y)$ and where the boundary density is $5$ and the boundary temperature follows a Gaussian distribution along the border with a standard deviation of $\sigma_b=0.05$:
\begin{equation}
    \forall z\in[0,1],~T_b(z)=0.5+9.5 e^{-\frac{z^2}{2\sigma_b}}.
\end{equation}
This boundary condition allows to model a beam of particles with a diameter of $\frac{\sqrt{2}}{2}$ at a high peak temperature ($T_b=10$) directed to the upper right corner. This test case is performed on structured meshes with the first and second order UGKS-M1 and with the HLL scheme (see figures \ref{fig:rayon_transv1} and \ref{fig:rayon_transv2}). On unstructured meshes, the first order UGKS-M1 solution is computed on a regular mesh and on an adapted mesh (see figure \ref{fig:uns}). From a physical point of view, we can notice that this stream of particle is transported through the low density area with little interactions in the transverse direction. In the high density area, the collisions predominate and induce energy deposition. The beam also causes a reduction in density through ablation along the front. This test case showcases a change of regime inside of the domain.

\begin{figure}
    \centering\input{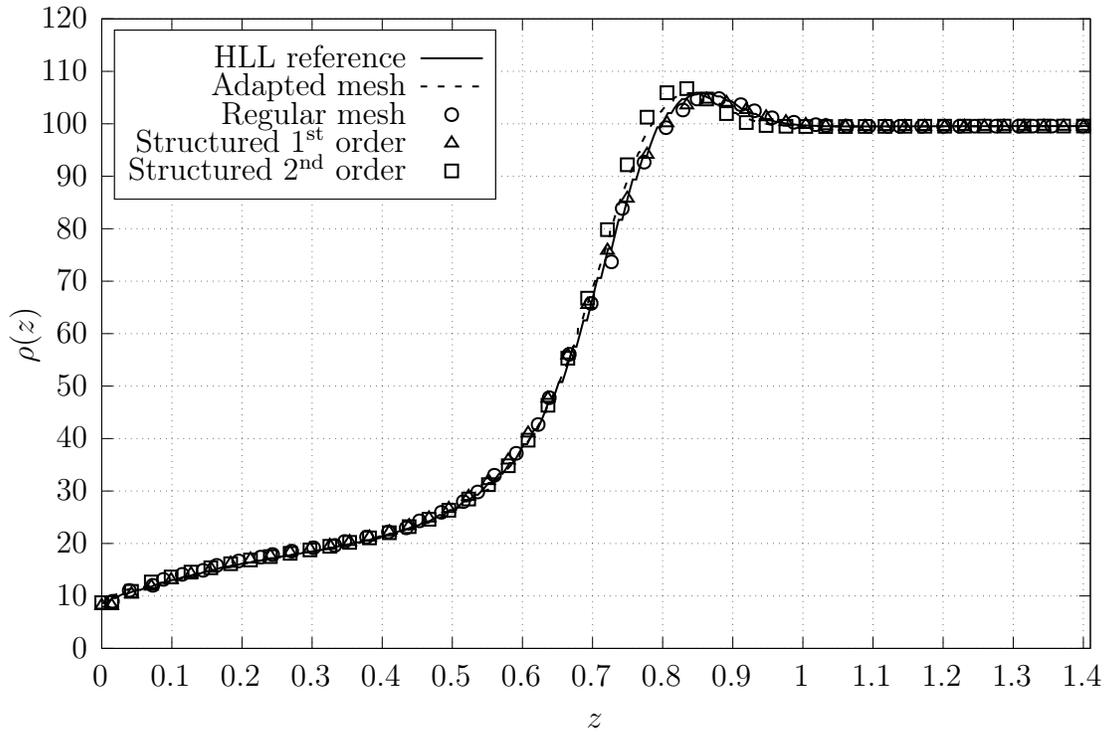}
    \input{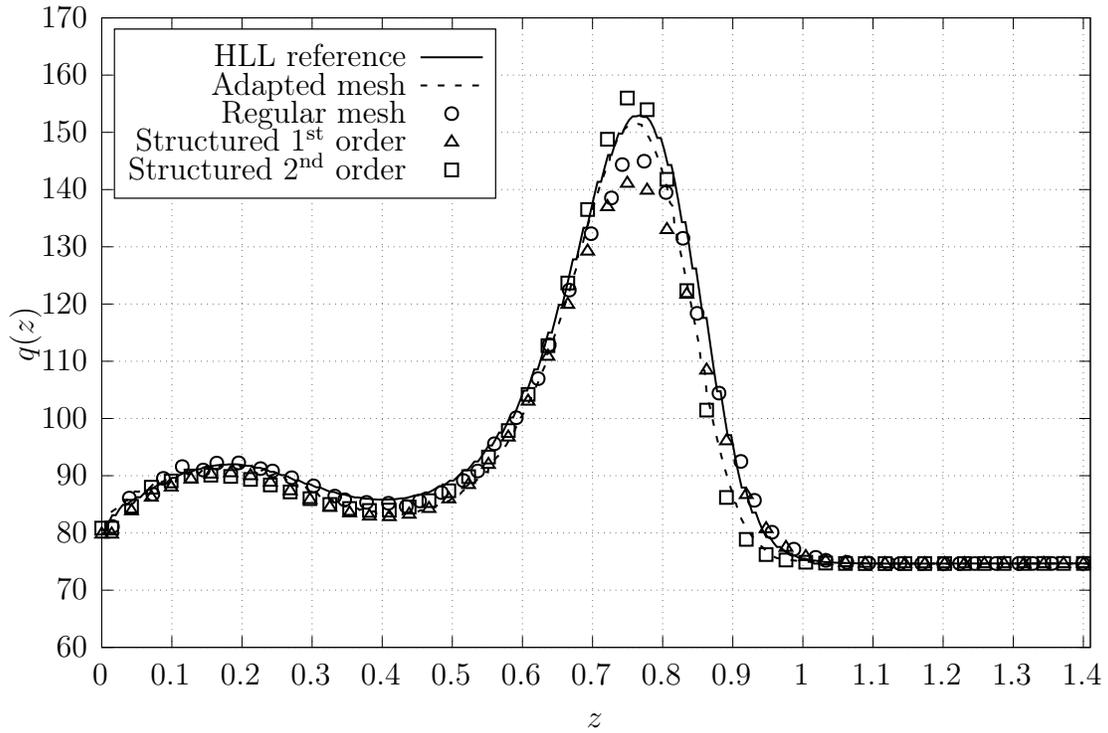}
    \caption{Test n°5: Energy deposition in a high density region. Density (above) and energy (below) in the transverse cutting line $y=x$ for the first order UGKS-M1 on structured and unstructured meshes (regular and adapted), the second order UGKS-M1 on structured mesh and for a first order HLL. \label{fig:uns_coupe}}
\end{figure}

From a numerical point of view, figures \ref{fig:rayon_transv1} and \ref{fig:rayon_transv2} show the UGKS-M1 solutions on structured meshes are similar to the one given by HLL but with less numerical diffusion. Moreover, the second order UGKS-M1 scheme allows to achieve the same accuracy than HLL but with 4 times less points in both directions. On unstructured meshes (see figure \ref{fig:uns}), the adapted mesh seems to significantly improve the solution as the energy maximum is greater than on the regular mesh. The areas of interest are overall well resolved. In order to make more accurate comparisons, both the density and the energy are represented on the cut line $y=x$ in figure \ref{fig:uns_coupe}. All the schemes seem to converge to the same solution on both the density and the energy. The numerical diffusion of each scheme is particularly visible on the energy maximum. The adapted mesh allows to significantly improve the solution quality at a reasonable cost.

\subsection{Test n°6: Non-local thermal transport}
The aim of this test case is to demonstrate the scheme capabilities in capturing non-local thermal transport effects in the context of inertial confinement fusion codes. In hydrodynamic codes, the electron heat flux needs to be closed. In order to do that, the stationary solution of the kinetic solution (\ref{equation_base}), in which the relaxation occurs towards an imposed hydrodynamic profile $\mf{\hat{W}}=\begin{pmatrix}\hat{\rho} & \hat{q}\end{pmatrix}^T$, can be computed at each hydrodynamic time step to evaluate the heat flux. From a numerical point of view, this can be achieved by setting the macroscopic fluxes to zero to prevent the evolution of the hydrodynamic variables. The stationary solution provides a non-local heat flux defined from the energy conservation law:
\begin{equation*}
    \partial_t q + \dfrac{1}{\eta} \boldsymbol\nabla_\mf{x} \cdot \ccrochet{v^3\boldsymbol\Omega f} = 0.
\end{equation*}
Thus in terms of the M1 variable, this flux is
\begin{equation}
    \boldsymbol\phi_q^{NL}=\dfrac{1}{\eta}\crochet{v^3 \mathbf{f}_1}.
\end{equation}
The non-local heat flux can be compared to the local one which is associated with the diffusion flux:
\begin{equation}
    \boldsymbol\phi_q^{L}=-\dfrac{10}{9\sigma} \boldsymbol\nabla_\textbf{x} \left(\dfrac{\hat{q}^2}{\hat{\rho}}\right)=-\dfrac{5}{2\sigma}\boldsymbol\nabla_\textbf{x} \left(\hat{\rho} \hat{T}^2\right).
\end{equation}
Due to its definition, the local heat flux only depends on the fixed hydrodynamic variables. The non-local flux takes into account for kinetic effects which induces non-local effects. 

\begin{figure}
    \centering\input{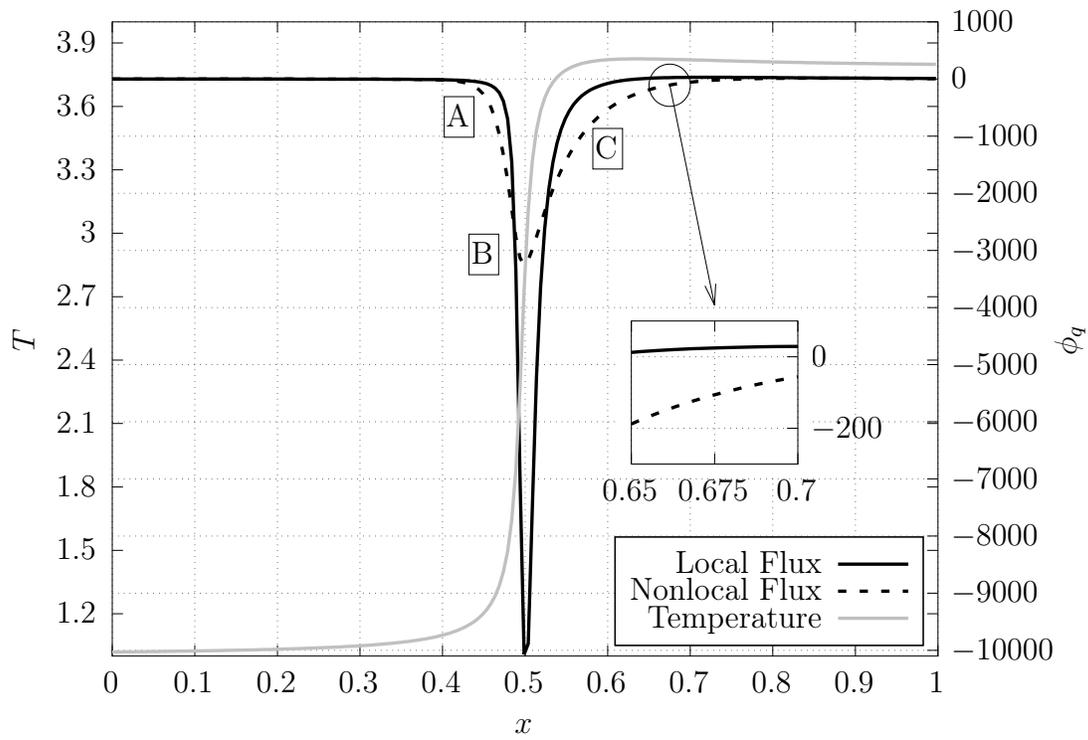}
    \caption{Test n°6: Nonlocal thermal transport. Temperature, local and nonlocal heat flux in the 1D domain given by UGKS-M1 in the mono-dimensional case. The letters highlights the main effects of nonlocal thermal transport: \textbf{A}-preheating, \textbf{B}-flux limitation, \textbf{C}-anti-natural flux and \textbf{D}-flux rotation.\label{fig:nl1}}
\end{figure}

In this test case, these non-local thermal transport effects are highlighted in 1D and 2D. First in the mono-dimensional case, we consider the relaxation of a regularised temperature step with a decrease on the right side. In figure \ref{fig:nl1}, the local and non-local heat fluxes are compared for the first order UGKS-M1 on a structured mesh. In front of the temperature step (\textbf{A}), a preheating effect can be notice as the non-local flux is non-zero despite the absence of temperature gradient. Moreover, in the strong gradient area (\textbf{B}), the non-local heat flux is lower than the local one by a factor of three. Finally, behind the temperature step (\textbf{C}), the local flux is positive since the temperature is slightly decreasing. However the non-local one is negative and as a consequence anti-natural. This kinetic effect is a consequence of the strong temperature gradient which directly influences the flux in this area.

\begin{figure}
    \centering\input{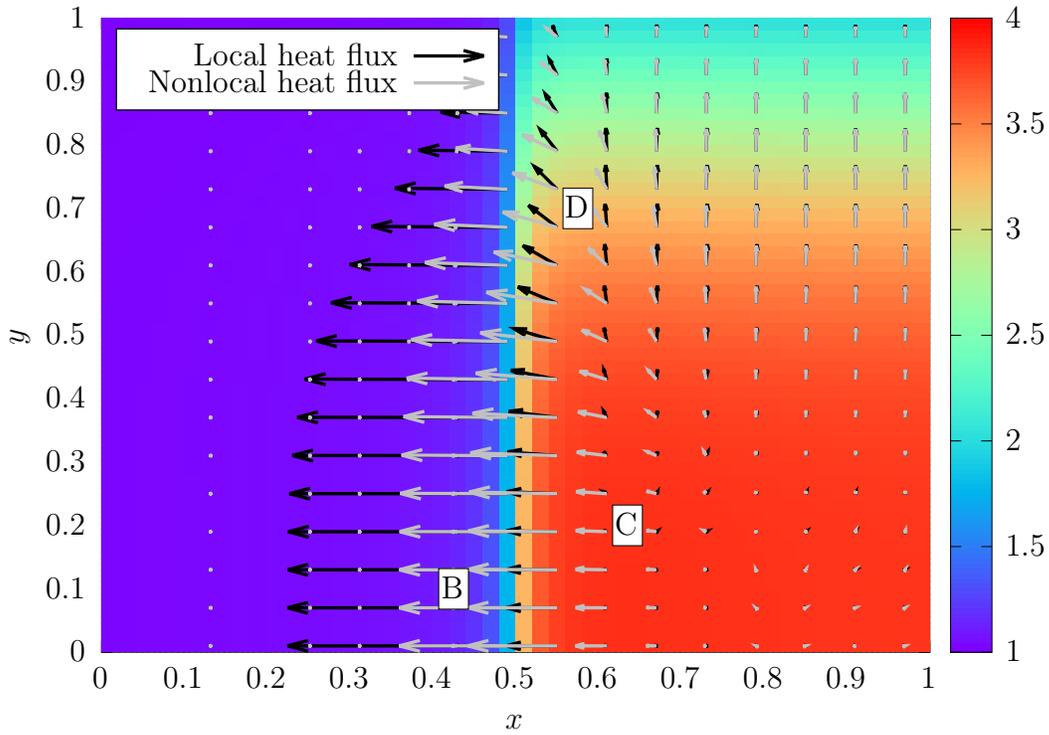}
    \caption{Test n°6: Nonlocal thermal transport. Temperature, local and nonlocal heat flux in the 2D domain given by UGKS-M1 in the bi-dimensional case. \label{fig:nl2}}
\end{figure}

In the bi-dimensional case, the same 1D profile is considered along the x-direction with a decay in the y-direction. The corresponding gradient is much smaller than the horizontal one to avoid inducing non-local effects, which explains why both fluxes coincide far from the temperature step in Figure \ref{fig:nl2}. On the contrary close to the step, the same 1D effects can be observed, especially the flux limitation (\textbf{B}) and the anti-natural flux (\textbf{C}). A flux rotation can also be observed on the right side of the step (\textbf{D}). The non-local flux is affected by the horizontal gradient which results in a rotation of the local flux towards this gradient despite the local gradient being in the vertical direction.

\section{Conclusion}
In this article, an asymptotic-preserving numerical scheme based on the Unified Gas Kinetic Scheme has been proposed for the M1 moment model of the electron transport equation. The method used consists in performing a moment closure in direction at the discrete level in the UGKS flux by using the M1 distribution function. It has been demonstrated that this procedure allows the scheme to inherit the asymptotic-preserving property of the UGKS and hence recovers correct numerical fluxes in the diffusion limit for the moment model scheme. As the fluxes are non-analytical due to the directional closure, a cost-efficient quadrature method that ensures the consistency of the scheme has been proposed. Moreover, a second-order extension that does not compromise the AP property has been suggested. It has also been demonstrated that several diffusion schemes can be easily recovered in UGKS and thus in the moment model scheme. Next, both schemes have been extended to unstructured meshes. The UGKS construction has been adapted in order to recover proper diffusion schemes in the limit. Two different versions of UGKS based on the diamond scheme have been proposed. Finally, several test cases have been chosen to validate the scheme in 1D and 2D, in every regime and on structured and unstructured meshes. The ability of the schemes to capture non-local thermal transport effects have also been demonstrated. This new scheme has also been compared with a standard HLL scheme and proven to be more accurate.

This article demonstrates that the UGKS can be used to obtain a good asymptotic-preserving scheme for a relevant moment model and proposes a general procedure to construct other asymptotic-preserving schemes. A natural extension of this work would be to improve the scheme on unstructured meshes by modifying the UGKS reconstructions in order to obtain other schemes in both the transport and diffusion regimes. Another perspective would be to extend this work to other moment models based on other kinds of collision kernels such as non-linear ones.

\bibliography{biblio}
\clearpage

\begin{figure}
    \centering\input{cinetique_rho} \\
    \centering\input{cinetique_T}
    \caption{Test n°2: Kinetic regime. Density and temperature in the domain at different times for the first order UGKS-M1 on a structured mesh, a regular and a deformed unstructured mesh along two line cuts C1 and C2. \label{fig:1dkinetic}}
\end{figure}
\begin{figure}
    \centering\input{intermediaire_rho}
    \centering\input{intermediaire_T}
    \caption{Test n°3: Intermediate regime. Density and temperature in the domain at different times for the first order UGKS-M1 on a structured mesh, a regular and a deformed unstructured mesh along two line cuts C1 and C2. \label{fig:1dinter2}}
\end{figure}
\begin{figure}
    \centering\input{diffusion_naif_rho}
    \centering \input{diffusion_naif_T}
    \caption{Test n°4: Diffusion regime. Density and temperature in the domain at different times for the first order UGKS-M1 with the naive diffusion scheme, on a structured mesh, a regular and a deformed unstructured mesh along two line cuts C1 and C2.\label{fig:diffusion1}}
\end{figure}
\begin{figure}
    \centering\input{diffusion_rho}
    \centering \input{diffusion_T}
    \caption{Test n°4: Diffusion regime. Density and temperature in the domain at different times for the first order UGKS-M1 with the diamond diffusion scheme, on a structured mesh, a regular and a deformed unstructured mesh along two line cuts C1 and C2.\label{fig:diffusion2}}
\end{figure}

\begin{figure}
\centering\scalebox{0.55}{
    \begin{tabular} {c c c}
    \LARGE{Reference solution (HLL 400 points)} & \LARGE{2\textsuperscript{nd} order UGKS-M1} \\
    \input{density_hll_fine} & \input{density_ugks2}  \\
    \LARGE{1\textsuperscript{st} order UGKS-M1} & \LARGE{HLL} \\
    \input{density_ugks} & \input{density_hll} \\
    \end{tabular}}
\caption{Test n°5: Energy deposition in a high density region. Density in the domain for the first and second order UGKS-M1 and HLL on structured meshes.\label{fig:rayon_transv1}}
\end{figure}
\begin{figure}
\centering\scalebox{0.55}{
    \begin{tabular} {c c c}
    \LARGE{Reference solution (HLL 400 points)} & \LARGE{2\textsuperscript{nd} order UGKS-M1} \\
    \input{energy_hll_fine} & \input{energy_ugks2} \\
    \LARGE{1\textsuperscript{st} order UGKS-M1} & \LARGE{HLL} \\
    \input{energy_ugks} & \input{energy_hll} \\
    \end{tabular}}
\caption{Test n°5: Energy deposition in a high density region. Energy in the domain for the first and second order UGKS-M1 and HLL on structured meshes.\label{fig:rayon_transv2}}
\end{figure}

\clearpage
\appendix

\section{M1 distribution moments \label{ann1}}
In this section, we recal the calculations performed in \cite{dubroca1999etude} in the context of radiative transfer. Let $\hat{f}(\boldsymbol\Omega)=e^{\alpha+\boldsymbol\beta\cdot\boldsymbol\Omega}$ be the M1 distribution function where $(\alpha,\boldsymbol\beta)\in\mathbb{R}\times\mathbb{R}^3$ are the entropic variables and $\boldsymbol\Omega \in \sphere$ is a vector on the sphere written in spherical coordinates:
\begin{equation*}
    \boldsymbol\Omega=\begin{pmatrix}
                           \sqrt{1-\mu^2}\cos \varphi \\
                           \sqrt{1-\mu^2}\sin \varphi \\
                           \mu
\end{pmatrix}, \text{ where } \left\{\begin{matrix}
\varphi \in [0,2\pi], \\
\mu \in [-1,1].
\end{matrix}\right.
\end{equation*}
Let $\mathbf{R}\in\mathscr{M}_3(\mathbb{R})$ be a rotation matrix such that $\mathbf{R}\boldsymbol\beta=||\boldsymbol\beta||\mathbf{e}_z$. The moments of the M1 distribution can be computed by performing the change of variable $\boldsymbol\Omega'=\mathbf{R}\boldsymbol\Omega$. The first moment is
\begin{equation}
\begin{aligned}
        f_0&=\int_{\sphere}e^{\alpha+\boldsymbol{\beta}\cdot\mathbf{\Omega}} \mathrm{d}\Omega, \\
        &=2\pi e^\alpha \int_{-1}^1 e^{||\boldsymbol\beta||\mu}\mathrm{d}\mu, \\
        &=4\pi e^\alpha \dfrac{\sinh{||\boldsymbol\beta||}}{||\boldsymbol\beta||}. \\
\end{aligned}
\end{equation}
This expression allows to get rid of the first entropic variable in the distribution function expression. The second moment is
\begin{equation}
\begin{aligned}
        \mathbf{f}_1&=\int_{\sphere}\mathbf{\Omega}e^{\alpha+\boldsymbol{\beta}\cdot\mathbf{\Omega}} \mathrm{d}\Omega, \\
        &=2\pi e^\alpha \int_{-1}^1 \mu e^{||\boldsymbol\beta|| \mu}\mathrm{d}\mu \mathbf{R}^T \mathbf{e}_z,\\
        &=4\pi e^\alpha \dfrac{\sinh{||\boldsymbol\beta||}}{||\boldsymbol\beta||} \left(\coth{||\boldsymbol\beta||}-\dfrac{1}{||\boldsymbol\beta||}\right) \dfrac{\boldsymbol\beta}{||\boldsymbol\beta||}, \\
        &=f_0 u \dfrac{\boldsymbol\beta}{||\boldsymbol\beta||},
    \end{aligned}
\end{equation}
where $u=\coth{||\boldsymbol\beta||}-\frac{1}{||\boldsymbol\beta||}$
This second expression allows to express the anisotropic variable $\boldsymbol\beta$ as a function of $\mf{U}=\begin{pmatrix}f_0 & \mf{f}_1\end{pmatrix}^T$. The third moment of the M1 distribution which provides the closure for the moment hierarchy is
\begin{equation}
    \begin{aligned}
        \mathbf{f}_2&=\int_{\sphere}\mathbf{\Omega}\otimes \mathbf{\Omega}e^{\alpha+\boldsymbol{\beta}\cdot\mathbf{\Omega}} \mathrm{d}\Omega, \\
        &=2\pi e^\alpha \int_{-1}^1 \mathbf{R}^T\left(\dfrac{1-\mu^2}{2} \mathbf{I}_3 + \dfrac{3\mu^2-1}{2}\mathbf{e}_z\otimes\mathbf{e}_z\right)\mathbf{R}e^{||\boldsymbol\beta|| \mu}\mathrm{d}\mu, \\
        &=f_0\dfrac{u}{||\boldsymbol\beta||}\mathbf{I}_3 + f_0\left(1-3\dfrac{u}{||\boldsymbol\beta||}\right) \dfrac{\boldsymbol\beta}{||\boldsymbol\beta||} \otimes \dfrac{\boldsymbol\beta}{||\boldsymbol\beta||}.
    \end{aligned}
\end{equation}

\section{M1 distribution half moments \label{ann2}}
Let $\beta_{xy}=\sqrt{\beta_x^2+\beta_y^2}$, $\mu_{xy}=\beta_{xy}\sqrt{1-\mu^2}$ and $I_n$ be the modified Bessel Function of the First Kind. The 1D expression of the half moments of the M1 distribution function are (for $k\geq 0$ and positive half moments):
\begin{subequations}
    \begin{align}
    \crochet{\Omega_z^k e^{\boldsymbol\beta \cdot \boldsymbol\Omega}\mathbb{1}_{\Omega_z\geq 0}}&=2\pi\int_{0}^1\mu^k e^{\beta_z \mu} \mathrm{I}_0(\mu_{xy})\mathrm{d}\mu, \\
   \crochet{\Omega_z^k \Omega_x\Omega_y e^{\boldsymbol\beta \cdot \boldsymbol\Omega}\mathbb{1}_{\Omega_z\geq 0}}&=\pi\int_{0}^1\mu^k(1-\mu^2) e^{\beta_z \mu} \sin{2\phi_0}\mathrm{I}_2(\mu_{xy})\mathrm{d}\mu, \\
    \crochet{\Omega_z^k\begin{pmatrix} \Omega_x \\ \Omega_y \end{pmatrix} e^{\boldsymbol\beta \cdot \boldsymbol\Omega}\mathbb{1}_{\Omega_z\geq 0}}&=2\pi\int_{0}^1\mu^k\sqrt{1-\mu^2} e^{\beta_z \mu} \begin{pmatrix} \cos{\phi_0} \\ \sin{\phi_0} \end{pmatrix}\mathrm{I}_1(\mu_{xy})\mathrm{d}\mu, \\
    \crochet{\Omega_z^k\begin{pmatrix} \Omega_x^2 \\ \Omega_y^2 \end{pmatrix} e^{\boldsymbol\beta \cdot \boldsymbol\Omega}\mathbb{1}_{\Omega_z\geq 0}}&=\pi\int_{0}^1\mu^k(1-\mu^2) e^{\beta_z \mu} \left(\mathrm{I}_0(\mu_{xy})+\begin{pmatrix} \cos{2\phi_0} \\ -\cos{2\phi_0} \end{pmatrix}\mathrm{I}_2(\mu_{xy})\right)\mathrm{d}\mu,
    \end{align}
\end{subequations}
where $\phi_0=\arctan{\dfrac{\beta_y}{\beta_x}}-(1-\sign{\beta_x})\dfrac{\pi}{2}$.

\section{Calculation of the macroscopic fluxes\label{ann0}}
The macroscopic fluxes $\boldsymbol{\Phi}_{i+1/2,j}$ are moments of the microscopic flux: $\ccrochet{\boldsymbol\Psi(v) \phi_{i+1/2,j}}$. For the calculation of these fluxes, it is necessary to express the moments of the slopes of the Maxwellian (in front of the term $D_{i+1/2,j}^n$ in (\ref{flux_micro})). This term is linked to the diffusion scheme since it remains the only term in the diffusion limit. First, the different moments of the Maxwellian are computed as functions of the conservative and entropic variables:
\begin{subequations}
\begin{align}
    \ccrochet{M[\mf{W}]}&=\rho, \\
    \ccrochet{v^2M[\mf{W}]}&=2q=3\rho T, \\
    \ccrochet{v^4M[\mf{W}]}&=5\dfrac{4}{3}\dfrac{q^2}{\rho}=15\rho T^2, \\
    \ccrochet{v^6M[\mf{W}]}&=35\dfrac{8}{9}\dfrac{q^3}{\rho^2}=105\rho T^3.
\end{align}
\end{subequations}
Next, these moments are utilized to compute the diffusion related term in each macroscopic fluxes. In the specific case where the Maxwellian is reconstructed by using the conservative variables $(\mf{W}_{i,j})$, the term in $\Phi^\rho_{i+1/2,j}$ is 
\begin{equation}
    \begin{aligned}
        D_{i+1/2,j}^n &\ccrochet{v^2\Omega_x^2\left(\delta_x^{Ln}M_{i+1/2,j}^n \mathbb{1}_{\Omega_x > 0}+\delta_x^{Rn}M_{i+1/2,j}^n \mathbb{1}_{\Omega_x < 0}\right)} \\
        &=\dfrac{D_{i+1/2,j}^n}{3\Delta x} {K}(\mf{W}_{i+1/2,j}^n)^T \ccrochet{v^2\boldsymbol\Psi(v) M[\mf{W}_{i+1/2,j}^n]} \cdot \left(\mf{W}_{i+1,j}^n-\mf{W}_{i,j}^n\right), \\
        &=\dfrac{2D_{i+1/2,j}^n}{3\Delta x}\begin{pmatrix} \dfrac{5}{2\rho_{i+1/2,j}^n} & -\dfrac{3}{2q_{i+1/2,j}^n} \\ -\dfrac{3}{2q_{i+1/2,j}^n} & \dfrac{3\rho_{i+1/2,j}^n}{2(q_{i+1/2,j}^n)^2}\end{pmatrix}\begin{pmatrix}q_{i+1/2,j}^n \\ \dfrac{5}{3}\dfrac{(q^n_{i+1/2,j})^2}{\rho_{i+1/2,j}^n}\end{pmatrix} \cdot \begin{pmatrix}\rho_{i+1,j}^n-\rho_{i,j}^n \\ q_{i+1,j}^n-q_{i,j}^n\end{pmatrix} ,\\
        &=\dfrac{2D_{i+1/2,j}^n}{3\Delta x}(q_{i+1,j}^n-q_{i,j}^n),
    \end{aligned}
\end{equation}
and the one in $\Phi^q_{i+1/2,j}$ is
\begin{equation}
    \begin{aligned}
        D_{i+1/2,j}^n &\ccrochet{\dfrac{1}{2}v^4\Omega_x^2\left(\delta_x^{Ln}M_{i+1/2,j}^n \mathbb{1}_{\Omega_x > 0}+\delta_x^{Rn}M_{i+1/2,j}^n \mathbb{1}_{\Omega_x < 0}\right)} \\
        &=\dfrac{D_{i+1/2,j}^n}{3\Delta x} {K}(\mf{W}_{i+1/2,j}^n)^T \ccrochet{\dfrac{1}{2}v^4\boldsymbol\Psi(v) M[\mf{W}_{i+1/2,j}^n]} \cdot \left(\mf{W}_{i+1,j}^n-\mf{W}_{i,j}^n\right), \\
        &=\dfrac{2D_{i+1/2,j}^n}{3\Delta x}\dfrac{5}{3}\begin{pmatrix} \dfrac{5}{2\rho_{i+1/2,j}^n} & -\dfrac{3}{2q_{i+1/2,j}^n} \\ -\dfrac{3}{2q_{i+1/2,j}^n} & \dfrac{3\rho_{i+1/2,j}^n}{2(q_{i+1/2,j}^n)^2}\end{pmatrix}\begin{pmatrix}\dfrac{(q^n_{i+1/2,j})^2}{\rho_{i+1/2,j}^n} \\ \dfrac{7}{3}\dfrac{(q^n_{i+1/2,j})^3}{(\rho_{i+1/2,j}^n)^2}\end{pmatrix} \cdot \begin{pmatrix}\rho_{i+1,j}^n-\rho_{i,j}^n \\ q_{i+1,j}^n-q_{i,j}^n\end{pmatrix},\\
        &=\dfrac{2D_{i+1/2,j}^n}{3\Delta x} \dfrac{5(q_{i+1/2,j}^n)^2}{3\rho_{i+1/2,j}^n}\left( -
    \dfrac{\rho_{i+1,j}^n-\rho_{i,j}^n}{\rho_{i+1/2,j}^n} + 2 \dfrac{q_{i+1,j}^n-q_{i,j}^n}{q_{i+1/2,j}^n}\right).
    \end{aligned}
\end{equation}
This procedure remains the same with a different Maxwellian reconstruction using another set of variables.

\section{Calculation of the Jacobian for the second order UGKS-M1 \label{ann4}}
First, the Jacobian of the inverse transformation $\mathbf{U} \to \boldsymbol\Lambda(\mathbf{U})$ is
\begin{equation}
    \mf{J}_\mf{U}(\boldsymbol\Lambda)=\begin{pmatrix} f_0 & \mf{f}_1^T \\ \mf{f}_1 & \mf{f}_2 \end{pmatrix}.
\end{equation}
To continue this calculation, the second moment of the M1 distribution function should be written in terms of the rotation matrix $\mf{R}$ involved in its calculation. A simple factorisation allows to show that
\begin{equation}
    \begin{aligned}
    \mf{J}_\mf{U}(\boldsymbol\Lambda)&=f_0\begin{pmatrix} 1 & \textbf{0} \\ \textbf{0} & \mathbf{R}^T \end{pmatrix}\begin{pmatrix} 1 & (\mathbf{R}\mathbf{u})^T \\ \mathbf{R}\mathbf{u} & \dfrac{u}{||\boldsymbol\beta||}\mathbf{I}_3 + (1-3\dfrac{u}{||\boldsymbol\beta||})\mathbf{e}_z \otimes \mathbf{e}_z\end{pmatrix} \begin{pmatrix} 1 & \textbf{0} \\ \textbf{0} & \mathbf{R} \end{pmatrix},\\
    &=       f_0\begin{pmatrix} 1 & \textbf{0} \\ \textbf{0} & \mathbf{R}^T \end{pmatrix}                       
                                    \begin{pmatrix} 1 & 0 & 0 & u \\ 
                                      0 & \dfrac{u}{||\boldsymbol\beta||} & 0 & 0 \\
                                      0 & 0 & \dfrac{u}{||\boldsymbol\beta||} & 0 \\
                                      u & 0 & 0 & 1-2\dfrac{u}{||\boldsymbol\beta||}\\
                                      \end{pmatrix}
                   \begin{pmatrix} 1 & \textbf{0} \\ \textbf{0} & \mathbf{R} \end{pmatrix}.                   
    \end{aligned}
\end{equation}
As the Jacobian of the inverse transformation is the inverse matrix, the Jacobian of $\boldsymbol\Lambda \to \mathbf{U}(\boldsymbol\Lambda)$ is
\begin{equation}
\begin{aligned}
    \mathbf{J}_\mathbf\Lambda (\mathbf{U})=\dfrac{f_0^{-1}}{\xi}
     \begin{pmatrix} 1-2\dfrac{u}{||\boldsymbol\beta||} & - \mathbf{u}^T \\ 
                                -\mathbf{u} & \dfrac{||\boldsymbol\beta||}{u}\xi \mathbf{I}_3+(1-\dfrac{||\boldsymbol\beta||}{u}\xi)\dfrac{\mathbf{u}}{||\mathbf{u}||}\otimes\dfrac{\mathbf{u}}{||\mathbf{u}||}\\
    \end{pmatrix},
    \end{aligned}
\end{equation}
where the determinant of the external matrix is $\xi=1-2\dfrac{u}{||\boldsymbol\beta||}-u^2$.

\section{Second order UGKS-M1 fluxes \label{ann5}}
The microscopic fluxes of the second order UGKS-M1 scheme are
\begin{subequations}
    \begin{align}
    \begin{split}
    \chi^0_{i+1/2,j}(v)=&A_{i+1/2,j}^n v \crochet{\Omega_x \hat{f}_{i,j}^{n(+)} \mathbb{1}_{\Omega_x\geq 0}+\Omega_x\hat{f}_{i+1,j}^{n(-)} \mathbb{1}_{\Omega_x\leq 0}} \\
     + & B_{i+1/2}^n v^2\crochet{ \Omega_x^2 {\delta_x}\hat{f}_{i,j}^n \mathbb{1} _{\Omega_x>0} + \Omega_x^2{\delta_x}\hat{f}_{i+1,j}^n \mathbb{1}_{\Omega_x<0}}\\
    +&\dfrac{D_{i+1/2,j}^n}{3\Delta x} v^2 \boldsymbol\Psi^T \mf{K}_0\left(\mf{W}_{i+1/2,j}^n\right) \left(\mf{W}_{i+1,j}^n-\mf{W}_{i,j}^n \right),
    \end{split} \\
    \begin{split}
    \boldsymbol\chi^1_{i+1/2,j}(v)=&A_{i+1/2,j}^n v \crochet{\Omega_x\boldsymbol\Omega \hat{f}_{i,j}^{n(+)} \mathbb{1}_{\Omega_x\geq 0}+\Omega_x\boldsymbol\Omega\hat{f}_{i+1,j}^{n(-)} \mathbb{1}_{\Omega_x\leq 0}} \\
     + & B_{i+1/2}^n v^2 \crochet{\Omega_x^2\boldsymbol\Omega {\delta_x}\hat{f}_{i,j}^n \mathbb{1} _{\Omega_x>0} +\Omega_x^2\boldsymbol\Omega{\delta_x}\hat{f}_{i+1,j}^n \mathbb{1}_{\Omega_x<0}}\\
    +&\dfrac{C_{i+1/2,j}^n}{3} v M[\mf{W}_{i+1/2,j}^n] \mathbf{e}_x \\
    -&\dfrac{D_{i+1/2,j}^n}{4\Delta x} v^2\boldsymbol\Psi^T \mf{K}_0\left(\mf{W}_{i+1/2,j}^n\right) \left(\mf{W}_{i+1,j}^n-2\mf{W}_{i+1/2,j}^n+\mf{W}_{i,j}^n \right)\mf{e}_x,
    \end{split}
    \end{align}
\end{subequations}
and the macroscopic fluxes are (when the conservative variables $\mf{W}$ are used to reconstruct the Maxwellian)
\begin{subequations}
    \begin{align}
    \begin{split}
    \Phi_{i+1/2,j}^\rho=&A_{i+1/2,j}^{n} \ccrochet{\indi{v\Omega_x \hat{f}_{i,j}^{n(+)}}{v\Omega_x\hat{f}_{i+1,n}^{n(-)}}{\Omega_x} } \\
    + & B_{i+1/2}^n \ccrochet{v^2 \Omega_x^2 {\delta_x}\hat{f}(\mathbf{U}_{i,j}^n) \mathbb{1} _{\Omega_x>0} + v^2\Omega_x^2{\delta_x}\hat{f}(\mathbf{U}_{i+1,j}^n) \mathbb{1}_{\Omega_x<0}}\\
    +&\dfrac{2D_{i+1/2,j}^n}{3\Delta x} (q_{i+1,j}^n-q_{i,j}^n), 
    \end{split} \\
    \begin{split}
    \Phi_{i+1/2,j}^q=&\dfrac{A_{i+1/2,j}^n}{2} \ccrochet{\indi{v^3\Omega_x \hat{f}_{i,j}^{n(+)}}{v^3\Omega_x \hat{f}_{i+1,n}^{n(-)}}{\Omega_x} } \\
    + & \dfrac{B_{i+1/2,j}^n}{2} \ccrochet{v^4 \Omega_x^2 {\delta_x}\hat{f}(\mathbf{U}_{i,j}^n) \mathbb{1} _{\Omega_x>0} + v^4\Omega_x^2{\delta_x}\hat{f}(\mathbf{U}_{i+1,j}^n) \mathbb{1}_{\Omega_x<0}}\\
    +&\dfrac{2D_{i+1/2,j}^n}{3\Delta x} \dfrac{5(q_{i+1/2,j}^n)^2}{3\rho_{i+1/2,j}^n}\left( -
    \dfrac{\rho_{i+1,j}^n-\rho_{i,j}^n}{\rho_{i+1/2,j}^n} + 2 \dfrac{q_{i+1,j}^n-q_{i,j}^n}{q_{i+1/2,j}^n}\right),
    \end{split}
    \end{align}
\end{subequations}
where $\hat{f}_{i,j}^{n(\pm)}=\hat{f}_{i,j}^n\pm \dfrac{\Delta x}{2}\delta_x\hat{f}_{i,j}^n$.

\section{Least square reconstruction \label{ann6}}
Following \cite{diamant} work, a least squared based method is used to reconstruct the macroscopic variables at the mesh cells vertices. Let $u$ be any macroscopic variable, $n$ be a mesh vertex and $\Tilde{u}_n(\mf{x})=a+\mf{b}\cdot(\mf{x}-\mf{x}_n)$ be an affine reconstruction of $u$ around $n$. The quadratic error $E=\frac{1}{2}\sum_{K\in\mathscr{T}_n} |\Tilde{u}_n(\mf{x}_{K})-u_{K}|^2$ is minimal if and only if the parameters of the reconstruction satisfies the following normal equations:
\begin{equation}
    \begin{pmatrix} p & \mathbf{I}^T \\ \mathbf{I} & \mathbf{J} \\ \end{pmatrix}
    \begin{pmatrix}a \\ \mf{b}\end{pmatrix}
    =\begin{pmatrix} \displaystyle\sum_{K\in\mathscr{T}_n} u_K \\ \displaystyle\sum_{K\in\mathscr{T}_n} (\mathbf{x}_K-\mathbf{x}_n)u_K \\ \end{pmatrix},
\end{equation}
where $\mathscr{T}_n$ is the set of cells that have $n$ as a vertex, $p\in \mathbb{N}$ is the cardinality of $\mathscr{T}_n$ and $(\mf{I},\mf{J})\in\mathbb{R}^2 \times \mathscr{M}_2(\mathbb{R})$ are a vector and a matrix defined as follows:
\begin{subequations}
\begin{align}
    \mathbf{I} &= \sum_{K\in\mathscr{T}_n}(\mathbf{x}_K-\mathbf{x}_n), \\
    \mathbf{J} &= \sum_{K\in\mathscr{T}_n}(\mathbf{x}_K-\mathbf{x}_n) \otimes (\mathbf{x}_K-\mathbf{x}_n). 
\end{align}
\end{subequations}
To determine the value of the reconstructed variable at the vertex $n$ (which is $\tilde{u}_n(\mf{x}_n)=a$), one must solve for $a$ in the previous equation. This quantity was derived by \cite{diamant} and can be written in the form:
\begin{equation}
    a=\sum_{K\in\mathscr{T}_n} \dfrac{1-\boldsymbol\lambda \cdot \mf{x}_K}{p-\boldsymbol\lambda\cdot\mf{I}} u_K,\quad \boldsymbol\lambda=\mf{J}^{-1}\mf{I}.
\end{equation}
From a numerical point of view, the reconstruction step at the mesh vertices consists in the computation of a dot product at each time step.

\end{document}